\documentclass[12pt,draftcls,onecolumn]{IEEEtran}

\usepackage{graphicx}
\usepackage{bm}
\usepackage{epsfig}
\usepackage{amssymb,amsfonts,amsmath}
\usepackage{subfigure}
\usepackage{setspace}
\usepackage{enumerate}
\usepackage{color}
\usepackage{dsfont}

\newtheorem{theorem}{Theorem}

\newtheorem{remark}{Remark}

\newtheorem{proposition}{Proposition}
\newtheorem{example}{Example}

\def\a {\alpha}
\def\b {\beta}
\def\e {\varepsilon}
\def\s {\sigma}
\def\w {\omega}
\def\r {\rho}
\def\g {\gamma}

\IEEEoverridecommandlockouts

\begin{document}

\title{\LARGE \bf
Ensemble Control of Time-Invariant Linear Systems with Linear Parameter Variation
}

\author{Jr-Shin~Li and Ji~Qi
\thanks{*This work was supported by the National Science Foundation under the award 1301148. The authors contributed equally to this work.}
\thanks{J.-S. Li and J. Qi are with the Department of Electrical and Systems Engineering, Washington University, St. Louis, MO 63130, USA
        (e-mail: {\tt\small jsli@ese.wustl.edu; qij@ese.wustl.edu).}}%
}


\maketitle

\begin{abstract}

In this paper, we study the control of a class of time-invariant linear ensemble systems whose natural dynamics are linear in the system parameter. This class of ensemble control systems arises from practical engineering and physical applications, such as transport of quantum particles and control of uncertain harmonic systems. We establish explicit algebraic criterions to examine controllability of such ensemble systems. Our derivation is based on the notion of polynomial approximation, where the elements of the reachable set of the ensemble system are represented in polynomials of the system parameter and used to approximate the desired state of interest. In addition, we highlight the role of the spectra of the system matrices play in the determination of ensemble controllability. Finally, illustrative examples and numerical simulations for optimal control of this class of linear ensemble systems are presented to demonstrate the theoretical results.
\end{abstract}


\section{Introduction}
Robust and sensorless manipulation of a collection of structurally similar dynamical systems with variation in common system parameters, or of a single system with uncertainty in the parameters, is compelling in various areas of science and engineering. Prominent examples range from the application of optimal pulses to produce a desired time evolution of a large quantum ensemble in quantum control \cite{Li_PRA06, Li_PNAS11, Glaser98}, and the use of external stimuli to desynchronize a population of neurons for the treatment of neurological disorders \cite{Benabid10, Foutz10, Li_JNE12}, to the implementation of open-loop controls for approximate steering of robots under bounded model perturbation in robotics \cite{Becker12rob}. Such practical control designs give rise to challenging problems involving the guidance of a large number or a continuum of structurally similar dynamical systems using a common open-loop control input, which arises because measurements for the state of each individual system of the ensemble is impractical and hence state feedback is unavailable.

The research in the control of ensemble systems has been active in both theoretical and computational aspects. The controllability for an ensemble of systems evolving on the Lie group SO(3) has been investigated through a conversion of the analysis to polynomial approximation \cite{Li_TAC09}. The necessary and sufficient controllability characterization of an ensemble of finite-dimensional time-varying linear systems was provided in terms of the singular system of the input-to-state operator that governs the system dynamics \cite{Li_TAC11}. Controllability and optimal control of an ensemble of weakly forced nonlinear oscillators, such as neuron and chemical oscillators, described by phase-reduced models were also analyzed \cite{Li_TAC13, Li_NOLCOS10, Li_JNE12, Li_PRL13, Wilson14}. Recently, a unified computational method for solving optimal ensemble control problems based on multidimensional pseudospectral approximations has been developed \cite{Li_TAC12_QCP} and successfully employed to design optimal pulses for protein NMR spectroscopy \cite{Li_PNAS11}. An optimization-free computational algorithm based on the singular value decomposition (SVD) was also established 
to compute minimum-energy controls for steering time-varying linear ensemble systems \cite{Li_ACC12SVD}, such as quantum transport systems \cite{Li_TAC14_Transport}. Numerous work in ensemble control has also emerged in the biological domain, with the aim to understand the coordination of the movement of flocks \cite{Brockett10flock} and to control large-scale complex networks \cite{Liu11, Hochberg06}.

Although intensive work has been conducted to characterize controllability of ensemble systems \cite{Li_TAC09, Beauchard10, Li_TAC11, Li_TAC13}, explicit and algebraically or numerically verifiable controllability conditions are still insufficiently explored. Our previous work illustrated that controllability of an ensemble of 
time-varying linear systems is determined by the growth rate of the singular values of the input-to-state Fredholm operator, 
which is, however, intractable to verify \cite{Li_TAC11}. In this paper, we study a class of time-invariant linear ensemble systems, whose natural dynamics are linear in the system parameter. We provide a detailed 
controllability characterization of such ensemble systems and derive explicit algebraic controllability conditions that are related to the system matrices. This class of linear ensemble systems arises from practical engineering and physical applications, such as the transport of an ensemble of atoms and the control of a harmonic oscillator with uncertain frequency \cite{Li_ACC12SVD, Li_TAC14_Transport}.

In the next section, we review existing fundamental results on ensemble control of linear systems, which motivate this work. In Section \ref{sec:main}, we construct the necessary and sufficient controllability conditions for the time-invariant linear ensemble systems of interest. Our derivation is based on the notion of polynomial approximation, where the elements of the reachable set are represented as polynomials of the system parameters and used to approximate the desired ensemble states. 
Examples and numerical simulations of optimal controls for steering such linear ensemble systems are illustrated in Section \ref{sec:example}.

\section{Fundamental Results on Ensemble Control of Linear Systems}
\label{sec:prelim}
In this section, we recapitulate key fundamental results on ensemble control of time-varying linear systems, which are pertinent to this study and motivate our theoretical developments. Meanwhile, through this review, we define the mathematical settings and notations that will be used throughout this paper. 

Consider an ensemble of finite-dimensional time-varying linear systems indexed by a parameter $\beta$ varying over a compact set $K\subset\mathbb{R}$, given by
\begin{eqnarray}
	\label{eq:linear_ensemble}
	\frac{d}{dt}{X}(t,\beta)=A(t,\beta)X(t,\beta)+B(t,\beta)u(t),
\end{eqnarray}
where $X\in M \subset\mathbb{R}^n$ is the state , $\beta\in K$, and $u\in L_2^m[0,T]$ is an $L_2$ control; the elements of $A(t,\beta)\in\mathbb{R}^{n\times n}$ and $B(t,\beta)\in\mathbb{R}^{n\times m}$ are complex $L_{\infty}$ and $L_2$ functions, respectively, defined on the compact set $D=[0,T]\times K$, and are denoted $A\in L_\infty^{n\times n}(D)$ and $B\in L_2^{n\times m}(D)$. We say that the system \eqref{eq:linear_ensemble} is \emph{uniformly ensemble controllable} on the underlying function space (in this case the space $L_2^n(K)$) if there exists a finite time $T>0$ and an open-loop control function $u:[0,T]\rightarrow\mathbb{R}^m$ (in this case an $L_2$ function) that steers the system from an initial state $X_0(\b)=X(0,\b)$ into an $\varepsilon$-neighborhood of a target state $X_F(\b)$ in time $T$, i.e., if $\sup_{\b\in K}\|X(T,\b)-X_F(\b)\|=\|X(T,\b)-X_F(\b)\|_{\infty}<\varepsilon$ holds for any $\varepsilon>0$ \cite{Li_TAC09}. Note that $T$ may depend on $\e$. 

\begin{remark}
	Ensemble controllability can also be defined according to the $L_p$-norms for $1\leq p\leq\infty$, namely, the system is $L_p$-ensemble controllable if
	$$\left(\int_{K}\|X(T,\b)-X_F(\b)\|^p \, d\b\right)^{\frac{1}{p}}<\e.$$
\end{remark}

The necessary and sufficient ensemble controllability conditions for the system \eqref{eq:linear_ensemble} in a Hilbert space setting have been derived, which are related to the solvability of the Fredholm integral operator that characterizes the system dynamics \cite{Li_TAC11}, 
given by 
\begin{align}
	\label{eq:L}
	(Lu)(\beta)=\int_0^{T}\Phi(0,\sigma,\beta)B(\sigma,\beta)u(\sigma)d\sigma=\xi(\beta),
\end{align}
where $\Phi(t,0,\beta)$ is the transition matrix for the homogeneous system $\dot{X}(t,\beta)=A(t,\beta)X(t,\beta)$ and $\xi(\beta)=\Phi(0,T,\beta)X_F(\beta)-X_0(\beta)$. The controllability conditions are represented in terms of the singular system of the operator $L$ as in \eqref{eq:L} and are given by \cite{Li_TAC11}
\begin{align}
	\label{eq:conditions}
	\text{(i)} \ \ \sum_{n=1}^{\infty}\frac{|\langle\xi,\nu_n\rangle_K|^2}{\sigma_n^2}<\infty,  \qquad \text{(ii)} \ \ \xi\in\overline{\mathcal{R}(L)},
\end{align}
where $(\sigma_n,\mu_n,\nu_n)$ is a singular system \cite{Gohberg03} of $L$ and $\overline{\mathcal{R}(L)}$ denotes the closure of the range space of $L$. 

The controllability characterization 
stated in \eqref{eq:conditions} is in terms of the growth rate of the singular values of $L$, 
and hence is intractable to verify even numerically although numerical calculations 
of the singular values and singular vectors can be efficient \cite{Li_ACC12SVD}. As a result, constructing controllability conditions that are practically checkable is compelling. 

\section{Controllability Conditions for Time-Invariant Linear Ensemble Systems}
\label{sec:main}
Consider the time-invariant linear ensemble system indexed by a parameter $\b$ varying on a compact set $K$, given by
\begin{eqnarray}
	\label{eq:LTI}
	\frac{d}{dt}{X}(t,\b)=A(\b)X(t,\b)+B(\b)U(t)=A(\b)X(t,\b)+\sum_{j=1}^m u_j(t) b_j(\b),	
\end{eqnarray}
where $X\in M \subset\mathbb{R}^n$, $\beta\in K\subset\mathbb{R^d}$, and $U:[0,T]\rightarrow\mathbb{R}^m$ is piecewise continuous; the matrices $A\in C^{n\times n}(K)$ and $B\in C^{n\times m}(K)$, whose elements are continuous functions over $K$, and $b_j$ is the $j^{th}$ column of $B$. Then, the input-to-state operator $\tilde{L}:PC^m[0,T]\rightarrow C^n(K)$ of the system \eqref{eq:LTI} is given by
\begin{align}
	\label{eq:L_1}
	(\tilde{L} u)(\beta)=\int_0^{T}\tilde{\Phi}(0,\sigma,\beta)B(\beta)U(\sigma)d\sigma=\tilde{\xi}(\beta),
\end{align}
where $\tilde{\Phi}(t,0,\beta)=e^{A(\b)t}$, 
$\tilde{\xi}(\beta)=\tilde{\Phi}(0,T,\b)X_F(\beta)-X_0(\beta)$,
and $PC^m$ and $C^n$ denote the space of $m$-tuples of piecewise continuous functions and $n$-tuples of continuous functions, respectively. It is known that the reachable set of the system \eqref{eq:LTI} starting from $X_0$ can be characterized by the range space of $\tilde{L}$, denoted as $R(\tilde{L})$, 
that is,
\begin{eqnarray}
	\label{eq:reachable}
	\mathcal{R}_T(X_0)=\tilde{\Phi}(T,0,\beta)(R(\tilde{L})+X_0).
\end{eqnarray}
In the following, we derive the relation between $R(\tilde{L})$ and the Lie algebra generated by the drift and control vector fields, 
which is essential to the characterization of ensemble controllability.

\begin{proposition}
	\label{prop:range}
	The closure of the range space of operator $\tilde{L}$ defined in \eqref{eq:L_1} coincides with that of the Lie algebra generated by the drift and control vector fields, which is given by
	\begin{eqnarray}
		\label{eq:L0}
		\mathcal{L}_0=\{A(\b)X,b_j\}_{LA}=\text{span}\left\{A^k(\beta)b_j,\, j=1,...,m;\,k=0,1,...\right\}.
	\end{eqnarray}
	That is, $\overline{R(\tilde{L})}=\overline{\mathcal{L}_0}$.
\end{proposition}
{\it Proof: } See Appendix \ref{appd:proof}. \hfill$\Box$ 

By Proposition \ref{prop:range}, we can express the reachable set in terms of $\mathcal{L}_0$ as
\begin{align}
	\label{reachbleset}
	\overline{\mathcal{R}_T(X_0)}=\tilde{\Phi}(T,0,\beta)(h+X_0), \quad h\in\overline{\mathcal{L}_0}.
\end{align}

\begin{proposition}
	\label{prop:ec}
	The system \eqref{eq:LTI} is uniformly ensemble controllable if and only if for any given initial, $X_0(\b)=X(0,\b)\in C^n(K)$, and target state, $X_F(\b)\in C^n(K)$, there exists a finite time $T>0$ such that
	$$\tilde{\xi}(\beta)=\tilde{\Phi}(0,T,\beta)X_F(\beta)-X_0(\beta)\in\overline{\mathcal{L}_0}.$$
\end{proposition}
{\it Proof: }
The systems \eqref{eq:LTI} is uniformly ensemble controllable, by definition and by \eqref{eq:reachable}, if and only if for any $X_0, X_F\in C^n(K)$ and any $\delta>0$, there exists some $h\in R(\tilde{L})$ and $T\in (0,\infty)$ such that
\begin{eqnarray}
	\label{eq:e1}
	\|X_F-\tilde{\Phi}(T,0,\beta)(h+X_0)\|_{\infty}<\delta.
\end{eqnarray}
Therefore,
\begin{align*}
	& \|\tilde{\Phi}(0,T,\b)X_F-X_0-h\|_{\infty} \\
	&= \|\tilde{\Phi}(0,T,\b)\big[X_F-\tilde{\Phi}(T,0,\b)(h+X_0)\big]\|_{\infty} \\
	&\leq \|\tilde{\Phi}(0,T,\b)\|_{\infty} \ \|X_F-\tilde{\Phi}(T,0,\b)(h+X_0)\|_{\infty} \\
	&< \delta \|\tilde{\Phi}(0,T,\b)\|_{\infty} \doteq \e,
\end{align*}
where the last inequality is due to \eqref{eq:e1} and the boundedness of $\|\tilde{\Phi}(t,0,\b)\|_{\infty}$, $\forall\,t\in [0,T]$ and $\forall\,\b\in K$.
We conclude that $\tilde{\xi}\in\overline{R(\tilde{L})}=\overline{\mathcal{L}_0}$. \hfill$\Box$

\begin{remark}
	The result of Proposition \ref{prop:ec} is evident since $\tilde{\xi}(\b)\in\overline{\mathcal{L}_0}=\overline{R(\tilde{L})}$ implies the existence of a solution $U\in PC^m [0,T]$ to the integral equation \eqref{eq:L_1}.
\end{remark}

\subsection{Controllability of Linear Ensemble Systems with Parameters across the Origin}
\label{sec:origin}
In this section, we study the class of time-invariant linear ensemble systems whose natural dynamics are linear in the system parameters, and construct explicit ensemble controllability conditions.

\begin{theorem} 
	\label{thm:ec_with_0}
	Consider the time-invariant linear ensemble system
	\begin{eqnarray}
		\label{eq:LTI_0}
		\Sigma_1: \left\{
		\begin{array}{ll}
			\frac{d}{dt}X(t,\beta)=\beta AX(t,\beta)+BU(t), \\
			\beta\in K=[-\b_1,\b_2]\subset\mathbb{R}, \quad \b_1, \b_2>0,
		\end{array}
		\right.
	\end{eqnarray}
	where $X\in M\subset\mathbb{R}^n$, 
	the control $U:[0,T]\rightarrow\mathbb{R}^m$ is a piecewise continuous function, and $A\in\mathbb{R}^{n\times n}$ and $B\in\mathbb{R}^{n\times m}$ are constant matrices. This system is uniformly ensemble controllable if and only if
	\begin{align}
		\label{eq:cond1}
		\text{(i)} & \quad\text{rank}(A)=n, \\
		\label{eq:cond2}
		\text{(ii)} & \quad\text{rank}(B)=n.
	\end{align}
	Note that condition (ii) implies that the number of control inputs $m$ is no less than the dimension of the system $n$.
	
\end{theorem}
{\it Proof:} We first rewrite the system \eqref{eq:LTI_0} as
$$\frac{d}{dt}X(t,\beta)=\beta AX(t,\beta)+\sum_{j=1}^m u_j b_j,$$
where $U=(u_1,\ldots,u_m)'$ and $b_j$ is the $j^{th}$ column of $B$. \vspace{-8pt}\\

\noindent (Sufficiency) Suppose that the conditions (i) and (ii) hold. The Lie algebra $\mathcal{L}_0$ defined as in \eqref{eq:L0} can be easily computed and is given by $\mathcal{L}_0=\text{span}\{\beta^kA^kb_j, \, j=1,\ldots,m, \, k=0,1,\ldots\}$.
It is then sufficient to show, according to Proposition \ref{prop:ec}, that for any given respective initial and target states, $X_0(\b)$ and $X_F(\b)$, the ensemble state $\xi(\b)=e^{-\b AT}X_F(\b)-X_0(\b)\in\overline{\mathcal{L}_0}$
for some $T\in (0,\infty)$. In other words, it is equivalent to showing that for any given $\e>0$, there exists an $\eta\in\mathcal{L}_0$ such that $\|\xi(\b)-\eta\|_{\infty}<\e$. Since $\eta\in\mathcal{L}_0$, it can be represented as a linear combination,
\begin{align}
	\eta &= \sum_{k=0}^{\infty} \sum_{j=1}^m \a_{jk}\b^k A^k b_j= \sum_{k=0}^{\infty}(\a_{1k}A^k b_1+\ldots+\a_{mk}A^k b_m)\b^k,
	\label{eq:eta}
\end{align}
where $\a_{jk}\in\mathbb{R}$ for $j=1,\ldots,m$ and $k=0,1,\ldots$. Also, because $\xi\in C^n(K)$, it can be uniformly approximated by a vector-valued polynomial of order $N$, i.e., $p_N(\b)=\sum_{k=0}^N c_k \b^k$, such that $\|\xi(\b)-p_N(\b)\|<\e$, where $c_k\in\mathbb{R}^n$ are coefficient vectors for $k=0,\ldots, N$. It remains to show that the coefficients $\a_{jk}$ as in \eqref{eq:eta} can be chosen so that $\a_{1k}A^k b_1+\ldots+\a_{mk}A^k b_m=c_k\doteq(c_{1k},c_{2k},\ldots,c_{nk})'$ for all $k=0,\ldots,N$. This is possible because $A$ and $B$ are of full rank, the underdetermined system of linear equations
$$\left[\begin{array}{c|c|c}A^k b_1 & \ldots & A^k b_m\end{array}\right]\left[\begin{array}{c}\a_{1k}\\ \vdots \\ \a_{mk}\end{array}\right]=\left[\begin{array}{c}c_{1k}\\ \vdots \\ c_{nk}\end{array}\right]$$
has a solution for all $k=0,1,\ldots,N$. \vspace{-8pt}\\

\noindent (Necessity) We will show that if either of the conditions in \eqref{eq:cond1} or \eqref{eq:cond2} fails to hold, then the system (\ref{eq:LTI_0}) is not ensemble controllable.

(Case I): Suppose that rank$(B)<n$ and that $B$ has a row of zeros, say, without loss of generality, the last row $\ell_n$. Thus, for the system with $\b=0$, the state that can be reached is of the form $\eta=\a_{j0} b_j\in\mathcal{L}_0$, $j=1,\ldots,m$, from \eqref{eq:eta}. 
Since $\ell_n=0$, the last entry of $\eta$ is zero. Therefore, the system \eqref{eq:LTI_0} is not ensemble controllable, because any given 
ensemble state $\xi(\b)$ with $\xi(0)$ having a nonzero last entry cannot be uniformly approximated by $\eta$. Alternatively, if $B$ has no rows of zeros, then we express the last row as a linear combination of the others, i.e., $\ell_n=\sum_{i=1}^{n-1}\alpha_i \ell_i$, where $\alpha_i\in\mathbb{R}$, $i=1,2,...n-1$, and at least one of them is nonzero. A simple row operation $T\in\mathbb{R}^{n\times n}$ applied to \eqref{eq:LTI_0}, where
\begin{eqnarray}
	\label{transT}
	T=\left[\begin{array}{c c c c c} 1 & 0 & 0 & \ldots & 0 \\ 0 & 1 & 0 & \ldots & 0 \\ \vdots \\ -\alpha_1 & -\alpha_2 & \ldots & -\alpha_{n-1} & 1 \end{array}\right],
\end{eqnarray}
results in a transformed system given by
\begin{eqnarray}
	\label{linsys3}
	\frac{d}{dt}(TX)=\beta(TAT^{-1})(TX)+(TB)U.
\end{eqnarray}
The last row of $TB$ contains all zeros, which is equivalent to the previous case.

(Case II): Suppose that rank$(A)<n$ and that $A$ has a row of zeros, say the last row. Then, the last row of the matrices $A^k B$, $k=1,2,\ldots$ contains only zeros. As a result, the last entry of any $\eta\in\mathcal{L}_0$ as in \eqref{eq:eta} is a constant function, and hence the system \eqref{eq:LTI_0} is not ensemble controllable. The case when $A$ has no rows of zeros can be shown in the same fashion as the corresponding case discussed in Case I. \hfill$\Box$

We now present several examples to demonstrate the rank conditions derived in Theorem \ref{thm:ec_with_0}.

\begin{example}
	\label{ex:harmonic}
	Consider steering a harmonic oscillator with uncertainty in its frequency $\omega$ from an initial state $X_0(\w)=X(0,\w)\in C^2(K)$ to a desired target state $X_F(\w)\in C^2(K)$, modeled by 
	\begin{eqnarray}
		\label{eq:harmonic}
		\frac{d}{dt}X(t,\omega)=A(\omega)X(t,\omega)+BU,
	\end{eqnarray}
	where
	$$A(\omega)=\omega A=\omega\left[\begin{array}{cc} 0 & -1 \\ 1 & 0 \\ \end{array}\right], \ B=\left[\begin{array}{cc} 1 & 0 \\0 & 1 \end{array}\right], \ U=\left[\begin{array}{c} u_1 \\ u_2 \end{array}\right],$$
	and the frequency is known to be in the range $\omega\in K=[-\nu,\nu]$ with $\nu>0$. 
	Since rank$(A)=2$ and rank$(B)=2$, this system is ensemble controllable according to the conditions \eqref{eq:cond1} and \eqref{eq:cond2}. This can be illustrated using the concept of polynomial approximation. The Lie algebra of the drift and control vector fields is given by
	\begin{eqnarray}
		\mathcal{L}_0=\text{span}\{\omega^k\left[\begin{array}{c} 1 \\0 \end{array}\right], \, \omega^k\left[\begin{array}{c} 0 \\1 \end{array}\right],\, k=0,1,2,\ldots\}.
	\end{eqnarray}
	It follows that for a desired $\xi(\omega)=(\xi_1(\w),\xi_2(\w))'=e^{-\omega AT}X_F(\omega)-X_0(\omega)$, there exists $P(\omega)=\sum_{k=1}^N c_k\omega^k$, where $c_k=(d_k,e_k)'\in\mathbb{R}^2$, such that $\|\xi(\w)-P(\w)\|_{\infty}<\e$, namely, $\xi_1(\w)\approx\sum_{k=1}^N d_k\w^k$ and $\xi_2(\w)\approx\sum_{k=1}^N e_k\w^k$ for $k=0,1,\ldots,N$.

If, however, there is only one control available, say $u_2=0$, then rank$(B)=1$ and
$$\mathcal{L}_0=\text{span}\{\omega^{2k}\left[\begin{array}{c} 1 \\0 \end{array}\right],\omega^{2k+1}\left[\begin{array}{c} 0 \\1 \end{array}\right],k=0,1,2,\ldots\}.$$
In this case, uniform approximation for a given $\xi(\w)=(\xi_1(\w),\xi_2(\w))'$ for $\w\in [-\nu,\nu]$ by an $\eta\in\mathcal{L}_0$, i.e., $\|\xi-\eta\|_{\infty}<\e$, is possible only when $\xi_1(\w)$ and $\xi_2(\w)$ are an even and an odd function, respectively. This verifies that the condition \eqref{eq:cond2} is necessary.
\end{example}

\begin{remark}
	\label{rmk:sufficient}
	The characterization of ensemble controllability stated in Theorem \ref{thm:ec_with_0} is for ensemble systems with parameter variation on a set across 
	zero. The situation becomes different if the zero parameter value is not included. It follows from the proof of Theorem \ref{thm:ec_with_0} that (i) and (ii) in \eqref{eq:cond1} and \eqref{eq:cond2} are sufficient controllability conditions, and (i) is a necessary condition for the ensemble system \eqref{eq:LTI_0} with $0\notin K$, i.e., $K\subset\mathbb{R}^+$ or $K\subset\mathbb{R}^-$. However, the condition (ii) is not necessary, which will be discussed in detail in Section \ref{sec:no_origin}.
\end{remark}

The following examples are made to illustrate the observation described in Remark \ref{rmk:sufficient}.

\begin{example}
	\label{ex:aircraft}
	Consider a linearized lateral-directional model that describes aircraft dynamics in the presence of uncertainty \cite{Lavretsky12},
	$$\frac{d}{dt}X(t,\epsilon)=A(\epsilon)X(t,\epsilon)+BU,$$
	where
	$$A(\epsilon)=\epsilon A=\epsilon\left[\begin{array}{ccc} 0 & 0.1 & -1 \\ 10 & 0.1 & 0 \\ 4 & 0 & 0.1 \end{array}\right], \quad B=\left[\begin{array}{ccc} 1 &  0 &  0 \\ 0 & 1 & 0 \\ 0 & 0 & 1 \end{array}\right].$$
	The state $X=(\gamma,p_s,r_s)'$, where $\gamma$ denotes the angle of sideslip, $p_s$ and $r_s$ denote the stability axis roll and yaw rate, respectively, and $U=(u_1,u_2,u_3)'$ is the control input. Unpredictable perturbations in the environment of the aircraft may result in dispersion in system dynamics, which we model $\epsilon\in K=[0.8,1.2]$. Because $\text{rank}(A)=3$ and $\text{rank}(B)=3$, these are sufficient for the system to be ensemble controllable. 
\end{example}

\begin{example}
	\label{ex:aircraft2}
	Recall the system in Example \ref{ex:aircraft} with different dynamics in $A$, given by
	$$A(\epsilon)=\epsilon\left[\begin{array}{ccc} 0 & 0.1 & 0 \\ 0 & 0 & 1 \\ 4 & 0 & 0 \end{array}\right],$$
	and with only one control $u_1$ available, i.e., $u_2=u_3=0$. In this case, $\text{rank}(A)=3$ but $\text{rank}(B)=1<3$. 
	However, this system is ensemble controllable. Observe that
	$$\mathcal{L}_0=\text{span}\Big\{\left[\begin{array}{c} \epsilon^{3k} \\ 0 \\ 0 \end{array}\right], \left[\begin{array}{c} 0 \\ 0\\ \epsilon^{3k+1} \end{array}\right], \left[\begin{array}{c} 0 \\ \epsilon^{3k+2} \\ 0 \end{array}\right]\Big\}$$
	for $k=0,1,2,\ldots$, then any given ensemble state $\xi\in C^3(K)$ over  $K=[0.8,1.2]\subset\mathbb{R}^+$ can be uniformly approximated by the polynomials synthesized using the vector fields in $\mathcal{L}_0$ according to the M\"{u}ntz$-$Sz\'{a}sz theorem (see Appendix \ref{appd:Muntz}). The fact of ensemble controllability indicates that the condition \eqref{eq:cond2} is not necessary when the parameter $0\notin K\subset\mathbb{R}^+$ (or $\mathbb{R}^-$). \hfill$\Box$
\end{example}

\subsection{Controllability of Linear Ensemble Systems with Positive or Negative Parameters}
\label{sec:no_origin}
Remark \ref{rmk:sufficient} along with Example \ref{ex:aircraft} and \ref{ex:aircraft2} motivates the need to reexamine the controllability conditions for the system \eqref{eq:LTI_0} when the values of the system parameter are strictly positive or negative, i.e., $0\notin K\subset\mathbb{R}^+$ (or $\mathbb{R}^-$). It is then essential to investigate the case for $\text{rank}(A)=n$ and $\text{rank}(B)=\nu<n$, and without loss of generality we will consider the case of $\text{rank}(B)=\nu=m$. 
Now recall the system \eqref{eq:LTI_0} and let $A=PJP^{-1}$ be the spectral decomposition of $A$, where $J$ is either a diagonal or a Jordan canonical matrix. Let $Y(t,\b)=P^{-1}X(t,\beta)$, then the system $Y(t,\b)$ follows $\frac{d}{dt}Y(t,\beta)=\beta JY(t,\beta)+P^{-1}BU(t)$. If $P^{-1}B$ has a row of zeros, say the $k^{th}$ row, then this system fails to be controllable because no controls act on the $k^{th}$ state variable. Hence, we focus on the case when $P^{-1}B$ has no rows of zeros. Since $\text{rank} (B)=m$, we may, without loss of generality, assume that the first $m$ rows of $P^{-1}B$ is an $m\times m$ invertible matrix, $B_0$. In addition, by a simple rearrangement of the term $(P^{-1}B)U$, we can put the system $Y$ into the form
\begin{eqnarray}
	\label{eq:Y}
	\frac{d}{dt}Y(t,\beta)=\beta JY(t,\beta)+\widetilde{B}V(t)
\end{eqnarray}
with the first $m$ rows of $\widetilde{B}$ form an $m\times m$ identity matrix 
and $V(t)=B_0 U(t)$.

In this section, we will derive explicit controllability conditions for the ensemble system \eqref{eq:LTI_0}, where we consider the cases when the system matrix $A$ is either diagonalizable or is similar to a Jordan matrix and the cases when $A$ has real or complex eigenvalues.

\subsubsection{$A$ is Diagonalizable} 
We first consider the case when $A\in\mathbb{R}^{n\times n}$ is diagonalizable and has real eigenvalues $\lambda_i$, $i=1,\ldots,n$. Hence, $J=\text{diag}(\lambda_1,\ldots,\lambda_n)$ is a real diagonal matrix. The following example serves as a motivating example that leads to the derivation of the controllability conditions.

\begin{example}[A Motivating Example]
	\label{ex:motivating}
	\rm Consider the ensemble system
\begin{eqnarray}
	\label{eq:Y1}
	\frac{d}{dt}Y(t,\beta)=\beta J_a Y(t,\beta)+\widetilde{B}V(t),
\end{eqnarray}
where $$J_a=\left[\begin{array}{ccc} 1 & 0 & 0 \\ 0 & 2 & 0 \\ 0 & 0 & a \end{array}\right],\quad \widetilde{B}=\left[\begin{array}{cc} 1 & 0 \\ 0 & 1 \\ 1 & 2 \end{array}\right],$$
$\beta\in K=[1,3]$, and $V(t)\in\mathbb{R}^2$.

(Case I: $a=0.5$) As shown in Proposition \ref{prop:ec}, this system is ensemble controllable if and only if for any given continuous ensemble state $\xi(\beta)=(\xi_1(\beta),\xi_2(\beta),\xi_3(\beta))'=e^{-\b J_aT}X_F-X_0$ on the interval $[1,3]$, it holds that
$\xi\in\overline{\mathcal{L}_0}$, where $\overline{\mathcal{L}_0}$ is the closure of the Lie algebra $\{\b J_aY,b_j\}_{LA}$ and $b_j$, $j=1,2$, is the $j^{th}$ column of $\widetilde{B}$. Because any element $\eta(\b)=(\eta_1(\b),\eta_2(\b),\eta_3(\b))'\in\mathcal{L}_0$ is of the form
\begin{eqnarray}
	\label{eq:ex4}
	\eta(\beta)=(\beta J)^k\sum_{k=0}^{\infty} \Big[c_k b_1+d_k b_2\Big]
	=\sum_{k=0}^{\infty} \Big\{c_k\left[\begin{array}{c} (\beta)^k \\ 0 \\ (0.5\beta)^k \end{array}\right]+d_{k}\left[\begin{array}{c} 0 \\ (2\beta)^{k} \\ 2(0.5\beta)^{k} \end{array}\right]\Big\},
\end{eqnarray}
where $c_k,d_{k}\in\mathbb{R}$, it is obvious that $\xi_1$ and $\xi_2$ can be uniformly approximated by $\eta_1$ and $\eta_2$, respectively, by appropriate choices of the coefficients $c_k$ and $d_{k}$. Therefore, what remains to be checked is whether the approximation $\|\xi_3(\b)-\eta_3(\b)\|_{\infty}<\e$ is possible whenever the coefficients $c_k$ and $d_k$ are chosen for approximating $\xi_1$ and $\xi_2$.

\begin{figure}[t]
	\centering
	\subfigure[]{\includegraphics[scale=0.55]{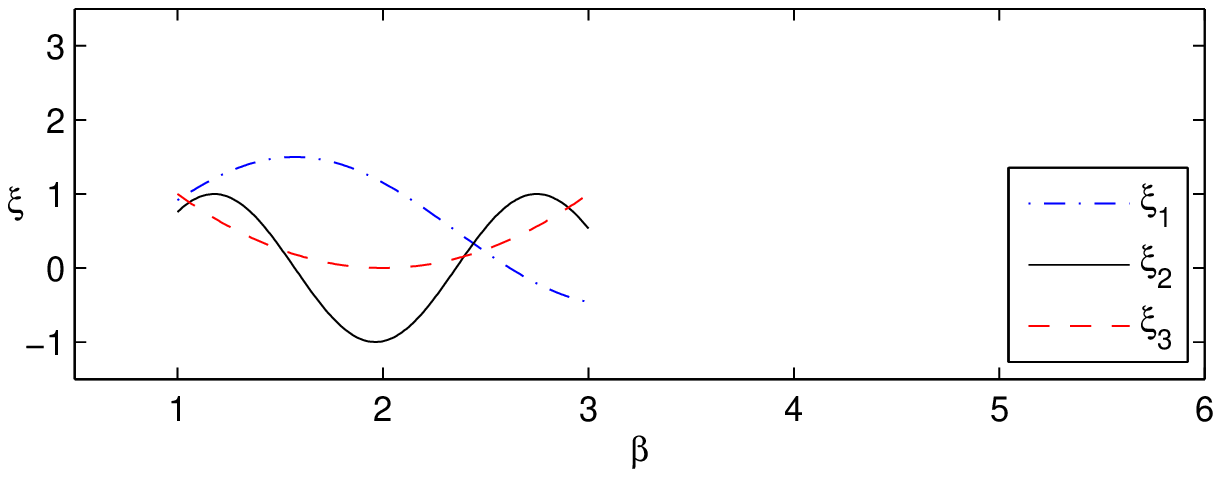}\label{fig:moti1}}
	\subfigure[]{\includegraphics[scale=0.55]{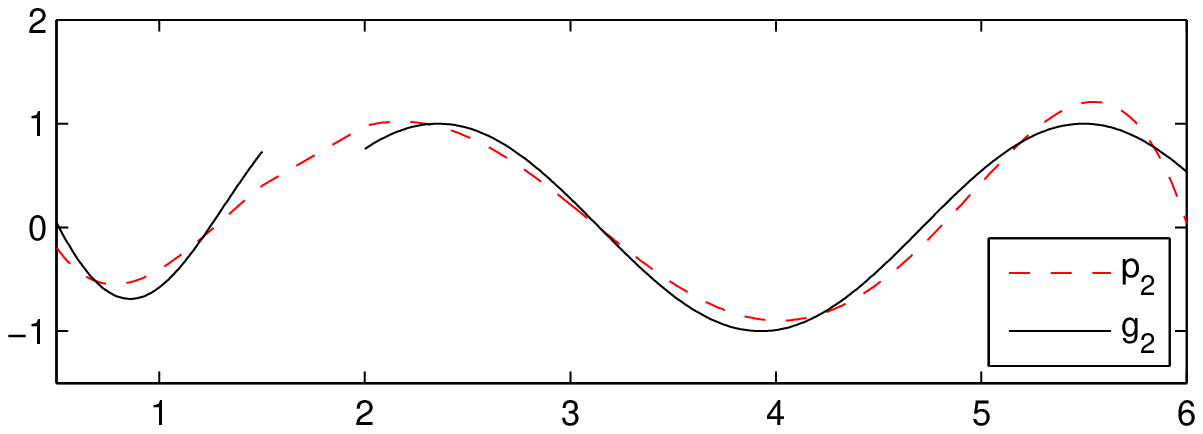}\label{fig:moti2}}
	\caption{An illustration of constructing an auxiliary continuous function based on a given desired ensemble state $\xi(\b)=(\xi_1(\b), \xi_2(\b), \xi_3(\b))'$. \subref{fig:moti1} The plot of $\xi_1(\beta)=-\cos(2\beta)+0.5$, $\xi_2(\beta)=-\sin(4\beta)$, and $\xi_3(\beta)=(\beta-2)^2$ for $\beta\in[1,3]$. \subref{fig:moti2} The constructed auxiliary function $g_2$ as in \eqref{eq:g2} approximated by a $5^{th}$ order polynomial $P_2=\sum_{k=0}^5 c_k\b^k$, where $c_0=3.27$, $c_1=-12.31$, $c_2=13.38$, $c_3=-5.68$, $c_4=1.02$, and $c_5=-0.07$.} 
     \label{fig:motivating}
\end{figure}

Since the eigenvalues of $J_a$ are $\lambda_1=1$, $\lambda_2=2$, and $\lambda_3=a=0.5$, the spectrum of $\b J_a$ is $\r(\b J_a)=\bigcup_{i=1}^3 \s_i$, where $\s_i=[\lambda_i,3\lambda_i]$. Let $\g_i=\lambda_i\b$ for $\b\in K$, then we have $\g_i\in\s_i$. Because $\s_1\cap \s_3=[1,1.5]\neq\emptyset$, $\eta_1(\b)$ and $\eta_3(\b)$ are correlated when $c_k$ are chosen. On the other hand, $\s_2\cap \s_3=\emptyset$, and thus $\eta_2=\sum_k d_k (2\b)^k=\sum_k d_k (\g_2)^k$ and $\eta_3=\sum_k [c_k(0.5\b)^k+ 2d_k(0.5\b)^k]=\sum_k (c_k+2d_k)(\g_3)^k$ in \eqref{eq:ex4} are two polynomials with different supports where $\g_2\in\s_2$ and $\g_3\in\s_3$. Therefore, it is possible to uniformly approximate $\xi_2(\b)$ and $\xi_3(\b)$, respectively, by a suitable choice of $d_k$ when $c_k$ are determined. This can be illustrated through the construction of the following auxiliary functions. Let $g_1(x)=\xi_1(x)$ for $x\in[1,3]$ and $g_2$ be a piecewise continuous function defined by
\begin{eqnarray}
	\label{eq:g2}
	g_2(x)=\left\{\begin{array}{ll}\frac{1}{2}\big[\xi_3(\frac{x}{0.5}) -\xi_1(\frac{x}{0.5})\big],  & x\in[0.5,1.5], \\ \xi_2(\frac{x}{2}), & x\in[2,6]. \end{array}\right.
\end{eqnarray}
Then, $g_1$ and $g_2$ can be uniformly approximated, respectively, by polynomials in $x$ such that $\|g_1(x)-\sum_k c_k x^k\|_{\infty}<\e_1$ for $x\in [1,3]$ and $\|g_2(x)-\sum_{k}d_{k} x^{k}\|_{\infty}<\e_2$ for $x\in [2,6]$. As a result, for any given function $\xi(\b)$ on the interval $[1,3]$, there exists an $\eta\in\mathcal{L}_0$ such that
$$\|\xi(\b)-\eta(\b)\|_{\infty}
=\Big\|\left[\begin{array}{c} \xi_1(\beta) \\ \xi_2(\beta) \\ \xi_3(\beta) \end{array}\right] - \left[\begin{array}{c} \sum_{k=0}^{\infty}c_k(\beta)^k \\ \sum_{k=0}^{\infty}d_{k}(2\beta)^{k} \\ \sum_{k=0}^{\infty}c_k(0.5\beta)^k+2\sum_{k=0}^{\infty}d_{k}(0.5\beta)^{k} \end{array}\right]\Big\|_{\infty}<\e,$$
and therefore the system (\ref{eq:Y1}) with $a=0.5$ is uniformly ensemble controllable. A visualization of this approximation is shown in Figure \ref{fig:motivating}.

(Case II: $a=1.5$) If $a=1.5$, then $\s_3=[1.5,4.5]$ and hence $\s_1\cap\s_2\cap\s_3=[2,3]\neq\emptyset$. In this case, no matter what $c_k$ and $d_{k}$ are chosen, $\xi_3$ is always a linear combination of $\xi_1$ and $\xi_2$, that is, $\xi_3(\beta)=\xi_1(\frac{3}{2}\beta)+2\xi_2(\frac{3}{4}\beta)$ for $\beta\in[\frac{4}{3},2]\subset [1,3]$. Hence, the system \eqref{eq:Y1} fails to be ensemble controllable. \hfill$\Box$
\end{example}

Example \ref{ex:motivating} highlights the role of the spectra of $A(\b)$ play in the determination of ensemble controllability and provides an insight into the development of uniform ensemble controllability conditions. Now consider the system \eqref{eq:Y}, the Lie algebra generated by the drift and control vector fields, which defines the reachable set of this ensemble system, is given by
$$\mathcal{L}_0=\text{span}\big\{\beta^kJ^k\tilde{b}_{j}, \, j=1,\ldots,m; \, k=0,1,2,\ldots \big\},$$
where $\tilde{b}_{j}$ is the $j^{th}$ column of $\widetilde{B}$. For any $\eta\in\mathcal{L}_0$, it is of the form
\begin{eqnarray}
	\label{eq:cmu}
	\left[\begin{array}{c}\eta_1 \\ \vdots \\ \eta_n\end{array}\right]
	=\sum_{k=0}^{\infty}\left[\begin{array}{ccc} \gamma_1^k & \cdots & 0 \\ \vdots & \ddots & \vdots \\  0 & \cdots & \gamma_n^k \end{array}\right]\widetilde{B}\left[\begin{array}{c} c_{1k} \\ \vdots \\ c_{mk} \end{array}\right],
\end{eqnarray}
where $\gamma_i=\lambda_i\b\in\s_i\doteq [\lambda_i\b_1,\lambda_i\b_2]$ for $i=1,\ldots, n$. Then, the spectra of $\b J$ for $\b\in K=[\b_1,\b_2]$ is given by $\rho(\b J)=\bigcup_{i=1}^n\s_i$. We further define $\hat{\eta}_i:\s_i\rightarrow\mathbb{R}$ by $\hat{\eta}_i(\g_i)=\eta_i(\frac{\g_i}{\lambda_i})=\eta_i(\b)$ 
for $i=1,\ldots,n$.

Let $\mathcal{C}=\{\mathcal{C}_j\}$ denote the collection of nonempty intersections of $\s_i$ for $i=1,\dots,n$. Furthermore, we define an equivalent set
\begin{eqnarray}
	\label{eq:setA}
	\mathcal{A}=\{\mathcal{A}_k\,|\,\bigcup_k\mathcal{A}_k=\bigcup_j\mathcal{C}_j\ \text{and}\ \bigcap_k\mathcal{A}_k=\emptyset\}.
\end{eqnarray}
In addition, let $\chi_k$ be the index set such that
\begin{eqnarray}
	\label{eq:chi}
	\chi_k=\{j\,|\, \mathcal{A}_k\subset\s_j,\, 1\leq j\leq n\}.
\end{eqnarray}
For instance, in Case I of Example \ref{ex:motivating}, we can choose $\mathcal{C}_1=\s_1\cap\s_3=[1,1.5]$ and $\mathcal{C}_2=\s_1\cap\s_2=[2,3]$. Since $\mathcal{C}_1$ and $\mathcal{C}_2$ are disjoint, we can choose $\mathcal{A}_1=\mathcal{C}_1$ and hence $\chi_1=\{1,3\}$, and $\mathcal{A}_2=\mathcal{C}_2$ and thus $\chi_2=\{1,2\}$. This gives $\mathcal{A}=\{\mathcal{A}_1,\mathcal{A}_2\}$.
Similarly, for Case II of Example \ref{ex:motivating}, we may choose $\mathcal{C}_1=\s_1\cap\s_2=[2,3]$,  $\mathcal{C}_2=\s_1\cap\s_3=[1.5,3]$, $\mathcal{C}_3=\s_2\cap\s_3=[2,4.5]$, and $\mathcal{C}_4=\s_1\cap\s_2\cap\s_3=[2,3]$. Then, we may pick $\mathcal{A}_1=\mathcal{C}_4=[2,3]$, $\mathcal{A}_2=\mathcal{C}_2\backslash\mathcal{A}_1=[1.5,2)$, and $\mathcal{A}_3=\mathcal{C}_3\backslash\mathcal{A}_1=(3,4.5]$ so that we have  $\bigcup_{k=1}^3\mathcal{A}_k=\bigcup_{j=1}^4\mathcal{C}_j=[1.5,4.5]$ and $\bigcap_{k=1}^3\mathcal{A}_k=\emptyset$, and obtain $\mathcal{A}=\{\mathcal{A}_1,\mathcal{A}_2,\mathcal{A}_3\}$, as well as $\chi_1=\{1,2,3\}$, $\chi_2=\{1,3\}$, and $\chi_3=\{2,3\}$.

Observe that if the functions $\hat{\eta}_j$, $j\in\chi_k$, are linearly dependent over $\mathcal{A}_k$, namely,
 \begin{eqnarray}
 	\label{eq:linrln}
 	\sum_{j\in\chi_k} a_j \hat{\eta}_j(\mu)=0,\quad \mu\in\mathcal{A}_k,
 \end{eqnarray}
for at least one $a_j\neq 0$, then the system \eqref{eq:Y} fails to be uniformly ensemble controllable, as shown in Case II of Example \ref{ex:motivating}. In addition, for any continuous function of the form $\xi(\b)=(\xi_1(\b),\ldots,\xi_m(\b),0,\ldots,0)'\in C^n(K)$ with $\xi_i\not\equiv 0$ for $i=1,\ldots,m$ and $0\notin K\subset\mathbb{R}^+$ (or $\mathbb{R}^-$), there exists an $\eta\in\mathcal{L}_0$ such that $\|\xi-\eta\|_{\infty}<\e$, because the subsystem $Y_m=(y_1,\ldots,y_m)'$ of $Y$ in \eqref{eq:Y}, which satisfies $\frac{d}{dt}Y_m(t,\b)=\b J_m Y_m(t,\b)+I_m V_m(t)$, where $J_m=\text{diag}(\lambda_1,\ldots,\lambda_m)$, $V_m=(v_1,\ldots,v_m)'$, and $I_m$ is the $m\times m$ identity matrix, is uniformly ensemble controllable according to Theorem \ref{thm:ec_with_0}. Consequently, efforts will be made to examine the reachability of the states $y_k$ for $k=m+1,\ldots,n$. 

To this end, we first introduce the collection of sets
\begin{eqnarray}
	\label{eq:Abar}
	\mathcal{\bar{A}}=\{A_r\subset\mathcal{A}\,|\, \max\{\chi_r\}>m\}.
\end{eqnarray}
For each set $\mathcal{A}_r\subset\mathcal{\bar{A}}$, i.e., the index set $\chi_r$ corresponding to $\mathcal{A}_r$ contains at least one index greater than $m$, we define an associated $n\times n$ binary, diagonal matrix,
\begin{eqnarray}
	\label{eq:Er}
	E_r=\text{diag}(\mathds{1}_1,\mathds{1}_2,\ldots,\mathds{1}_n),
\end{eqnarray}
where $\mathds{1}_i=1$ if $i\in\chi_r$ and $\mathds{1}_i=0$ if $i\notin\chi_r$ for $i=1,\ldots,n$. We further define the matrix $M_r=E_r\widetilde{B}\in\mathbb{R}^{n\times m}$, and then let $\overline{M}_r\in\mathbb{R}^{p_r\times m}$, $p_r\leq n$, be the matrix defined as $M_r$ with its rows of zeros removed if there are any. For example, in Case I of Example \ref{ex:motivating}, we have $n=3$, $m=2$, and $\bar{\mathcal{A}}=\{\mathcal{A}_1\}$ since $\chi_1=\{1,3\}$. This leads to
\begin{eqnarray}
	\label{eq:Mr}
	E_1=\left[\begin{array}{ccc}1 & 0 & 0 \\ 0 & 0 & 0 \\ 0 & 0 & 1\end{array}\right], \
M_1=\left[\begin{array}{cc}1 & 0 \\ 0 & 0 \\ 1 & 2\end{array}\right], \
\overline{M}_1=\left[\begin{array}{cc}1 & 0 \\ 1 & 2\end{array}\right],
\end{eqnarray}
and hence $p_1=2<n$. In the following, we show that the controllability of the system \eqref{eq:Y} is characterized by the rank of $\overline{M}_r$.

\begin{theorem}
	\label{thm:diagonal}
	Given the time-invariant linear ensemble system
	\begin{eqnarray}
		\label{eq:LTI_1}
		\Sigma_2:\left\{
		\begin{array}{ll}
			\frac{d}{dt}X(t,\beta)=\beta AX(t,\beta)+BU(t), \\
			\beta\in K=[\b_1,\b_2]\subset\mathbb{R}^+ \ (\text{or}\ \mathbb{R}^-),
		\end{array}
		\right.
	\end{eqnarray}
	where $X\in M\subset\mathbb{R}^n$, the control $U:[0,T]\rightarrow\mathbb{R}^m$ is piecewise continuous, and $A\in\mathbb{R}^{n\times n}$ and $B\in\mathbb{R}^{n\times m}$ are constant matrices with $\text{rank}(A)=n$ and $\text{rank}(B)=m<n$. Suppose that $A$ is diagonalizable and $A=PJP^{-1}$ is the eigen-decomposition with $J=\text{diag}(\lambda_1,\ldots,\lambda_n)$, where $\lambda_i\in\mathbb{R}$, $i=1,\ldots,n$. Consider the transformed system $Y(t,\b)=P^{-1}X(t,\b)$, which satisfies $\frac{d}{dt}Y(t,\beta)=\beta JY(t,\beta)+\widetilde{B}V(t)$, where $\widetilde{B}$ and $V(t)$ are defined in \eqref{eq:Y}. Let $S_k=\{j\,|\, \widetilde{B}_{kj}\neq 0, 1\leq j\leq m, m+1\leq k\leq n\}$ be an index set associated with the $k^{th}$ row of $\widetilde{B}$, which denotes that the $k^{th}$ state, $k\geq m+1$, receives the controls $u_j$, $j=1,\ldots,m$. The system \eqref{eq:LTI_1} is uniformly ensemble controllable if and only if
	\begin{enumerate}
		\item[(i)] $S_k\neq\emptyset$ for all $k=m+1,\ldots,n$, and
		\item[(ii)] $\bar{\mathcal{A}}=\emptyset$, or when $\bar{\mathcal{A}}\neq\emptyset$, for all $\mathcal{A}_r\subset\bar{\mathcal{A}}$,
			\begin{equation}
				\label{eq:conds2}
				\text{rank}(\overline{M}_r)=p_r\leq m,
			\end{equation}
			where $\bar{\mathcal{A}}$ and $\overline{M}_r\in\mathbb{R}^{p_r\times m}$ are defined above as in \eqref{eq:Abar} and \eqref{eq:Mr}. 
	\end{enumerate}
\end{theorem}

{\it Proof:} \noindent (Sufficiency)
Suppose that the conditions (i) and (ii) hold. We wish to show that any given ensemble state associated with a specified pair of initial and target states, i.e., $\xi=(\xi_1,\ldots,\xi_n)'=e^{-\b AT}X_F(\b)-X_0(\b)\in C^n (K)$, is reachable, namely, $\xi\in\overline{\mathcal{L}_0}$, where $\mathcal{L}_0=\{\b AX,b_i\}_{LA}=\{\b JX,\tilde{b}_i\}_{LA}$, and $b_i$ and $\tilde{b}_i$, $i=1,\ldots,m$, are the $i^{th}$ columns of $B$ and $\widetilde{B}$, respectively. Here, we show the case of $0<\b_1<\b_2$, and the proof of its counterpart, $\b_1<\b_2<0$, follows the same procedure.

Let $\s_i=\lambda_i K=[\lambda_i\b_1,\lambda_i\b_2]$ and let $\g_i=\lambda_i\b\in\s_i$. 
We first define a rescaling function, $\hat{\xi}$, on $\s_i$ such that $\hat{\xi}_i(\g_i)=\xi_i(\frac{\g_i}{\lambda_i})=\xi_i(\b)$. Similarly, for any $\eta=(\eta_1,\ldots,\eta_n)'\in\mathcal{L}_0$, we define $\hat{\eta}_i(\g_i)=\eta_i(\frac{\g_i}{\lambda_i})=\eta_i(\b)$ for $i=1,\ldots,n$. We now provide a constructive proof to show the sufficient condition. The procedure is based on constructing a sequence of auxiliary functions that are used to approximate the desired ensemble states. To begin with, we define the functions $g_l$ restricted on $\s_l$ such that $g_l|{\s_l}=\hat{\xi}_l|{\s_l}$ and set the domains $D_l=\s_l$ for $l=1,\ldots,m$. The ensemble controllability is then analyzed through the characterization of the set $\bar{\mathcal{A}}$, defined in \eqref{eq:Abar}, associated with the system parameters and the intersections of the spectra $\s_i$ for $i=1,\ldots,n$.

\noindent (Case I: $\mathcal{\bar{A}}=\emptyset$) In this case, the spectrum $\s_k$ for $k>m$ does not intersect with any other set $\s_i$, i.e., $\s_{k}\cap\s_i=\emptyset$ for $k=m+1,\ldots,n$, $i=1,\ldots,n$, and $k\neq i$. We start the proof with considering the $(m+1)^{th}$ row of $\widetilde{B}$. Let $k=m+1$ and since $S_k\neq\emptyset$, we pick $\a\in S_k$ (note $\a\leq m$) and extend the functions $g_l$ on $\s_k$ such that
\begin{eqnarray}
	\label{eq:g}
	\left\{
	\begin{array}{ll}
		g_j|\s_k=\frac{1}{\widetilde{B}_{kj}}(\hat{\xi}_k|\s_k), & j=\a \\
		g_j|\s_k=0, & j\in S_k,\quad j\neq\a,
	\end{array}
	\right.
\end{eqnarray}
where $\widetilde{B}_{kj}$ is the $kj^{th}$ element of $\widetilde{B}$. Note that now each function $g_j$, $j\in S_k$, is piecewise continuous on the extended domain $D_j^{(k)}=\s_j\cup\s_k$. By the Weierstrass approximation theorem, there exist $c_{lq}\in\mathbb{R}$ with $l=1,\ldots,m$ and $q\in\mathbb{N}$, such that
\begin{eqnarray}
	\label{eq:g_approxi}
	\|g_l(\g)-\sum_{q=0}^{N(\e)} c_{lq}\g^q\|_{\infty}<\e,
\end{eqnarray}
for $\g\in D_l$, where $N(\e)\in\mathbb{Z}^+$ depends on $\e$. This leads to the fact that for any given $\xi\in C^n(K)$, there exists an $\hat{\eta}\in\mathcal{L}_0$ such that $\|\hat{\eta}-\hat{\xi}\|_{\infty}<\e$, where $\hat{\eta}$ and $\hat{\xi}$ are defined above. Because when $1\leq l\leq m$, we have for $\g_l\in\s_l$,
$$\|\hat{\eta}_l(\g_l)-\hat{\xi}_l(\g_l)\|_{\infty}=\|\sum_{q=0}^{N(\e)}\sum_{j=1}^m c_{jq}\g_l^q\widetilde{B}_{lj}-\hat{\xi}_l(\g_l)\|_{\infty}=\|\sum_{q=0}^{N(\e)} c_{lq}\g_l^q-\hat{\xi}_l(\g_l)\|_{\infty}<\e,$$
by using \eqref{eq:g_approxi} with the definitions $\widetilde{B}_{lj}=\delta_{lj}$ for $1\leq j\leq m$ and $\hat{\xi}_l=g_l$ on $\s_l$, where $\delta_{lj}$ is the Kronecker delta; and for $k=m+1$, 
the same coefficients can be used to obtain 
\begin{align}
	\|\hat{\eta}_k(\g_k)-\hat{\xi}_k(\g_k)\|_{\infty}&=\|\sum_{q=0}^{N(\e)}\sum_{j=1}^m c_{jq}\g_k^q\widetilde{B}_{kj}-\hat{\xi}_k(\g_k)\|_{\infty} \nonumber
	\label{eq:k>m}
	=\|\sum_{q=0}^{N(\e)}\sum_{j\in S_k} c_{jq} \g_k^q\widetilde{B}_{kj}-\hat{\xi}_k(\g_k)\|_{\infty}\\
	&\leq\sum_{\substack{j\in S_k \\ j\neq\a}}\|\sum_{q=0}^{N(\e)} c_{jq} \g_k^q\widetilde{B}_{kj}\|_{\infty}+\|\sum_{q=0}^{N(\e)} c_{\alpha q} \g_k^q\widetilde{B}_{k\a }-\hat{\xi}_k(\g_k)\|_{\infty} \nonumber \\
	& < \e\sum_{j\in S_k}\widetilde{B}_{kj}.
\end{align}
We now repeat the same procedure presented above for $k=m+2$ to $k=n$. In each step $k$, the functions $g_j$, for $j\in S_k$, are extended to the corresponding $\s_k$ as described in \eqref{eq:g} and defined over the extended domain $D_j^{(k+1)}=D_j^{(k)}\cup\s_k$. This concludes that $\hat{\xi}\in\overline{\mathcal{L}_0}$, and hence, by Proposition \ref{prop:ec}, the system is uniformly ensemble controllable.

\noindent (Case II: $\mathcal{\bar{A}}\neq\emptyset$) For each $\mathcal{A}_r\subset\mathcal{\bar{A}}$ with $\chi_r=\{\ell_1,\ldots,\ell_j\}$, $j\leq m$, and the associated $E_r$ as defined in \eqref{eq:chi} and \eqref{eq:Er}, respectively, we have $\text{rank}(E_r)=|\chi_r|$, the cardinality of the index set $\chi_r$, and thus $\overline{M}_r\in\mathbb{R}^{|\chi_r|\times m}$ with $\text{rank}(\overline{M}_r)=|\chi_r|=p_r\leq m$ since \eqref{eq:conds2} holds. This follows that the subsystem of the system $Y$ in \eqref{eq:Y},
\begin{equation*}
	\frac{d}{dt}{\left[\begin{array}{c}y_{\ell_1} \\ y_{\ell_2} \\ \vdots \\ y_{\ell_j}\end{array}\right]}=\left[\begin{array}{cccc} \gamma_{\ell_1} & 0 & \cdots & 0 \\ 0 & \gamma_{\ell_2} &  \cdots & 0 \\ \vdots & \quad\quad\ddots \\  0 & \cdots & 0 & \gamma_{\ell_j} \end{array}\right]\left[\begin{array}{c}y_{\ell_1} \\ y_{\ell_2} \\ \vdots \\ y_{\ell_j}\end{array}\right]+\overline{M}_rU,
\end{equation*}
is ensemble controllable over the set $\mathcal{A}_r$ according to Remark \ref{rmk:sufficient}. Now let $\widetilde{\s}_k=\s_k\backslash\cup_{\mathcal{A}_r\in\mathcal{\bar{A}}}\mathcal{A}_r$ for $k=m+1,\ldots,n$. Because $\widetilde{\s}_k$'s are disjoint, so are the sets $\mathcal{A}_r\subset\bar{\mathcal{A}}$, we conclude, based on the result of Case I, that the system \eqref{eq:LTI_1} is uniformly ensemble controllable when the condition \eqref{eq:conds2} holds.

\noindent (Necessity) We will prove the negation of the condition by showing that when either (i) or (ii) is violated, then the system \eqref{eq:LTI_1} is not ensemble controllable. If the condition (i) does not hold, then the system \eqref{eq:LTI_1} is not ensemble controllable, because it implies that there exists at least one state $y_k$, where $m+1\leq k\leq n$, that receives no control inputs. If the condition (ii) is not satisfied, namely, there exists some $\mathcal{A}_r\in\bar{\mathcal{A}}$ with $\text{rank}(\overline{M}_r)<p_r$, then for any $\mu\in\mathcal{A}_r$, the rows of $\overline{M}_r$ are linearly dependent, namely, $q_{\ell_k}(\mu)=\sum_{\ell_i\in\chi_r, \ell_i\neq\ell_k} a_i q_{\ell_i}(\mu)$, where $q_{\ell_i}$ is the $\ell_i^{th}$ row of $\overline{M}_r$ and $a_i\in\mathbb{R}$ with at least one $a_i\neq 0$. Observe that the same linear dependence exists among the elements of $\hat{\eta}\in\mathcal{L}_0$ through the relation $\hat{\eta}_{\ell_k}(\mu)=\sum_{\ell_i\in\chi_r, \ell_i\neq \ell_k} a_i\hat{\eta}_{\ell_i}(\mu)$ for $\mu\in\mathcal{A}_r$. This implies that an arbitrary ensemble state $\hat{\xi}$ may not be reached, which leads to the failure of ensemble controllability.  \hfill$\Box$

\begin{remark}[Single-Input Systems]
	\label{rmk:single_input}
If the ensemble system \eqref{eq:LTI_1} receives a single input, i.e., $U:[0,T]\rightarrow\mathbb{R}$, then, from Theorem \ref{thm:diagonal}, $\mathcal{\bar{A}}=\emptyset$ if \eqref{eq:LTI_1} is uniformly ensemble controllable. This implies that there are no repeated eigenvalues among $\beta A$, or equivalently $\beta J$, where $\b\in K=[\b_1,\b_2]\subset\mathbb{R}$, which coincides with the previous result stated in \cite{Li_TAC09}. 
\end{remark}

\begin{remark}[Complexity for Verifying the Controllability Condition]
	\label{rmk:complexity}
	In Theorem \ref{thm:diagonal}, because each spectrum $\s_i$, $i=1,\ldots,n$, is a connected set, the upper bound of the cardinality of $\mathcal{A}$, i.e., the maximum possible number of the disjoint sets $\mathcal{A}_k$ defined in \eqref{eq:setA}, is $2n-1$. In other words, the complexity of examining the ensemble controllability condition (ii) in Theorem \ref{thm:diagonal} for the system \eqref{eq:LTI_1} is linear in the system dimension.	
\end{remark}

Now, we use the following example to illustrate the rank condition provided in Theorem \ref{thm:diagonal}.

\begin{example} 	
	\label{ex:thm2}
	Consider the ensemble system \eqref{eq:LTI_1} with
	$$A=\left[\begin{array}{cccc} 1 & 0 & 0 & 0 \\ 0 & 6 & 0 & 0 \\ 0 & 0 & \alpha & 0 \\ 0 & 0 & 0 & 2.5 \end{array}\right],\quad B=\left[\begin{array}{cc} 1 & 0 \\ 0 & 1 \\ 1 & 2 \\ 1 & 0 \end{array}\right],$$
	and $\b\in K=[1,2]$. We thus have the spectra $\s_i=\lambda_i K$ of $\b A$, which are given by $\s_1=[1,2]$, $\s_2=[6,12]$, $\s_3=[\alpha,2\alpha]$ and $\s_4=[2.5,5]$, where $\lambda_i$, $i=1,\ldots,4$, are the eigenvalues of $A$. We consider two values of $\a$ to demonstrate different scenarios characterized by $\bar{\mathcal{A}}$ as described in the proof of Theorem \ref{thm:diagonal}.

\noindent (Case I: $\alpha=0.4$) In this case, $\s_3=[0.4,0.8]$, and note that $n=4$ (four states) and $m=2$ (two controls). We then can define the index sets $S_3=\{1,2\}\neq\emptyset$ and $S_4=\{1\}\neq\emptyset$, and obtain $\mathcal{\bar{A}}=\emptyset$. By Theorem \ref{thm:diagonal}, the system is uniformly ensemble controllable. We now construct auxiliary functions following the procedures described in the proof of Theorem \ref{thm:diagonal} to illustrate the controllability through polynomial approximation. Consider steering the system between an initial and a target state, $X_0(\b)=X(0,\b)$ and $X_F(\b)$, respectively, in time $T$. The associated ensemble state is then obtained,    given by $\xi(\b)\doteq (\xi_1,\xi_2,\xi_3,\xi_4)'=\Phi(0,T,\b)X_F(\b)-X_0(\b)\in C^4(K)$, where $\Phi(t,0,\b)$ is the transition matrix related to the system \eqref{eq:LTI_1}. Let $\hat{\xi}=(\hat{\xi}_1,\hat{\xi}_2,\hat{\xi}_3,\hat{\xi}_4)$ be the rescaling of $\xi(\b)$ such that $\hat{\xi}_i:\s_i\rightarrow\mathbb{R}$ and $\hat{\xi}_i(\gamma_i)=\xi_i(\frac{\gamma_i}{\lambda_i})=\xi_i(\b)$ for $i=1,\ldots,4$. Because there are two controls ($m=2$), we define two functions $g_1$ and $g_2$ restricted on $\s_1$ and $\s_2$, respectively, by $g_1|\s_1=\hat{\xi}_1|\s_1$ and $g_2|\s_2=\hat{\xi}_2|\s_2$, and denote $D_1=\s_1$ and $D_2=\s_2$. Next, for $k=m+1=3$, we extend $g_1$ and $g_2$ onto $\s_3$ by letting $g_1|\s_3=\hat{\xi}_3|\s_3$ and $g_2|\s_3=0$, and set $D_1^{(3)}=\s_1\cup\s_3$ and $D_2^{(3)}=\s_2\cup\s_3$. Finally, for $k=n=4$, we extend $g_1$ onto $\s_4$ such that $g_1|\s_4=\hat{\xi}_4|\s_4$ and set $D_1^{(4)}=D_1^{(3)}\cup\s_4$. Since $g_2$ is not extended to $\s_4$ due to $S_4=\{1\}$, we have $D_2^{(4)}=D_2^{(3)}$. After such expansions, $g_1$ and $g_2$ are now piecewise continuous on $D_1^{(4)}$ and $D_2^{(4)}$, respectively, given by
\begin{eqnarray*}
	g_1(x)=\left\{\begin{array}{ll}\hat{\xi}_3(x),  & x\in[0.4,0.8], \\ \hat{\xi}_1(x), & x\in[1,2], \end{array}\right. \quad
	g_2(x)=\left\{\begin{array}{ll} 0,  & x\in[0.4,0.8], \\ \hat{\xi}_4(x) & x\in[2.5,5], \\ \hat{\xi}_2(x), & x\in[6,12]. \end{array}\right.
\end{eqnarray*}
It is then easy to verify that given any $\xi\in C^4(K)$, these functions $g_1$ and $g_2$ can be synthesized and then approximated by an $\eta\in\mathcal{L}_0$ 
of the form as in \eqref{eq:cmu}, so that $\|\hat{\eta}-\hat{\xi}\|_{\infty}<\e$ and therefore $\|\eta-\xi\|_{\infty}<\e$, where $\hat{\eta}_i(\g_i)=\eta_i(\frac{\g_i}{\lambda_i})=\eta_i(\b)$ with $\g_i\in\s_i$ and $\b\in K$. \vspace{4pt}

\noindent (Case II: $\alpha=4$) In this case, there are two sets of nonempty intersections between $\s_i$ and $\s_k$ for $i=1,\ldots,4$ and $k=3,4$, namely, $\mathcal{\bar{A}}=\{\mathcal{A}_1, \mathcal{A}_2\}$, where $\mathcal{A}_1=\s_3\cap\s_4=[4,5]$ and $\mathcal{A}_2=\s_2\cap\s_3=[6,8]$, and we denote the respective index sets $\chi_1=\{3,4\}$ and $\chi_2=\{2,3\}$. We also obtain $$\overline{M}_1=\left[\begin{array}{cc} 1 & 2 \\ 1 & 0 \end{array}\right]\quad\text{and} \quad \overline{M}_2=\left[\begin{array}{cc} 0 & 1 \\ 1 & 2 \end{array}\right].$$
Since $\overline{M}_1$ and $\overline{M}_2$ satisfy the rank condition \eqref{eq:conds2}, and $S_3=\{1,2\}$ and $S_4=\{1\}$ are both nonempty, this system is uniformly ensemble controllable. Similarly, we illustrate a construction of auxiliary functions that are used to demonstrate the ensemble controllability. Initially, following the same step as in Case I, we define two functions $g_1$ and $g_2$ restricted on $\s_1$ and $\s_2$, respectively, by $g_1|\s_1=\hat{\xi}_1|\s_1$ and $g_2|\s_2=\hat{\xi}_2|\s_2$, and denote $D_1=\s_1$ and $D_2=\s_2$. Then, we extend $g_1$ and $g_2$ onto $\mathcal{A}_1$ by letting $g_1|\mathcal{A}_1=\hat{\xi}_4|\mathcal{A}_1$ and $g_2|\mathcal{A}_1=\frac{1}{2}(\hat{\xi}_3|\mathcal{A}_1-\hat{\xi}_1|\mathcal{A}_1)$, and set $D_1^{(2)}=\s_1\cup\mathcal{A}_1$ and $D_2^{(2)}=\s_2\cup\mathcal{A}_1$. Furthermore, we extend $g_1$ onto $\mathcal{A}_2$ by $g_1|\mathcal{A}_2=\hat{\xi}_3|\mathcal{A}_2-2\hat{\xi}_2|\mathcal{A}_2$ and set $D_1^{(3)}=D_1^{(2)}\cup\mathcal{A}_2$ (note that $D_2^{(3)}=D_2^{(2)}$). Next, we extend, for $k=m+1=3$, $g_1$ and $g_2$ onto $\widetilde{\s}_3=\s_3\backslash (\mathcal{A}_1\cup\mathcal{A}_2)=(5,6)$ through a continuous function $f$ by letting $g_1|\widetilde{\s}_3=f|\widetilde{\s}_3$ and $g_2|\widetilde{\s}_3=\frac{1}{2}(\hat{\xi}_3|\widetilde{\s}_3-f|\widetilde{\s}_3)$, where $f$ is chosen so that $g_1$ is continuous on the closure of $\widetilde{\s}_3$, and set $D_1^{(4)}=D_1^{(3)}\cup\widetilde{\s}_3$ and $D_2^{(4)}=D_2^{(3)}\cup\widetilde{\s}_3$. Finally, for $k=n=4$, we extend $g_1$ to $\widetilde{\s}_4$ by letting $g_1|\widetilde{\s}_4=\hat{\xi}_4|\widetilde{\s}_4$ and $D_1^{(5)}=D_1^{(4)}\cup\widetilde{\s}_4$, where $\widetilde{\s}_4=\s_4\backslash\mathcal{A}_1=[2.5,4)$ (note that $D_2^{(5)}=D_2^{(4)}$). The constructed piecewise continuous functions $g_1$ and $g_2$ can be approximated with arbitrary accuracy by an $\eta\in\mathcal{L}_0$, as described in Case I, and consequently any given ensemble state $\xi\in C^4(K)$ can be approximately reached. \hfill$\Box$
\end{example}

\subsubsection{$A$ is in a Jordan Form}
\label{sec:jordan}
Here, we consider the case when $A$ in \eqref{eq:LTI_1} is similar to a Jordan form, given by $J=\text{diag}(J_{1},J_{2},...,J_{k})$, where $J_i=\text{diag}(J_{i1},\ldots,J_{ir})$, $i=1,\ldots,k$, and $J_{ij}$, $j=1,\ldots,r(i)$, denotes a Jordan block with real eigenvalues $\lambda_i\in\mathbb{R}$. Note that $J_{ij}$ can be a scalar. We first derive the ensemble controllability conditions for the case when $J$ is a pure Jordan block.

\begin{proposition}
	\label{prop:jordan}
	Given the time-invariant linear ensemble system
	\begin{eqnarray}
		\label{eq:jordanprop}
		\frac{d}{dt}X(t,\beta)=\beta J_0 X(t,\beta)+BU(t),
	\end{eqnarray}
	where
	$X\in\mathbb{R}^n$, $\b\in K=[\b_1,\b_2]\subset\mathbb{R}^+$ (or $\mathbb{R}^-$), $U:[0,T]\rightarrow\mathbb{R}^m$ is piecewise continuous, 
	$$J_0=\left[\begin{array}{ccccc} \lambda & 1 &  & & \\  & \lambda & 1 & & \\  &  & \ddots & \ddots &  \\  & & & \ddots & 1\\ &  &  & & \lambda \end{array}\right]$$
	with $\lambda\in\mathbb{R}$, and $B\in\mathbb{R}^{n\times m}$ is a constant matrix. This system is uniformly ensemble controllable if and only if $\text{rank}(B)=n$.
\end{proposition}

{\it Proof:}
The sufficiency of the rank condition directly follows the proof to the sufficiency of Theorem \ref{thm:ec_with_0}. We will prove the necessity by showing that if $\text{rank}(B)<n$, then the system (\ref{eq:jordanprop}) fails to be ensemble controllable. It is sufficient to focus on the case when $\text{rank}(B)=n-1$ and $B\in\mathbb{R}^{n\times(n-1)}$.

We first illustrate the result through the case of $n=2$, where we consider $B=(1,\a)'$ with $\a\neq 0$. In this case, we have
$$J_0(\lambda)=\left[\begin{array}{cc} \lambda & 1 \\ 0 & \lambda \\ \end{array}\right]\quad\text{and} \quad J_0^k(\lambda)=\left[\begin{array}{cc} \lambda^k & k\lambda^{k-1} \\ 0 & \lambda^k \\ \end{array}\right].$$
The reachable set of the system is characterized by $\mathcal{L}_0=\{\b J_0 X,b_j\}_{LA}$, where $b_j$, $j=1,\ldots,m$, is the $j^{th}$ column of $B$, and hence the states that can be 
reached are $\eta=(\eta_1,\eta_2)'\in\mathcal{L}_0$ and of the form
$$\eta_1(\b)=\sum_{k=0}^{N} c_k(\beta\lambda)^k+\a\sum_{k=1}^{N} c_k k\beta^k\lambda^{k-1},\quad \eta_2(\b)=\a\sum_{k=0}^{N} c_k(\beta\lambda)^k.$$
Observe that for any $N>0$ as long as the coefficients $c_k$, $k=0,1,2,\ldots$, are determined, the states of the system that can be approximated are correlated through the relation $\eta_1=\frac{1}{\a}\eta_2+\frac{d}{d\lambda}\eta_2$. In other words, given a desired ensemble state $\xi\in C^2(K)$ and a specified $\e>0$, there may exist no $\{c_k\}$ and $N$ such that $\|\eta-\xi\|_{\infty}<\e$, which concludes that the system \eqref{eq:jordanprop} is not ensemble controllable.

For $n>2$, we will show that $\eta_1$ is correlated to $\eta_2,\ldots,\eta_n$ and their higher order derivatives, and hence the system \eqref{eq:jordanprop} is not ensemble controllable. We assume, without loss of generality, that the first $n-1$ rows of $B$ form an identity matrix and suppose that
$$B=\left[\begin{array}{ccccc} 1 \\  &  1 \\  &  & \ddots \\ & & & 1 \\ a_1 & a_2 & \cdots & a_{n-1}\end{array}\right]\in\mathbb{R}^{n\times (n-1)},$$
where $a_1\neq 0$. Note that if $a_1=0$, it is equivalent to studying the system \eqref{eq:jordanprop} with a reduced dimension. The states of the system \eqref{eq:jordanprop} that can be reached, i.e., $\eta\in\mathcal{L}_0$, are of the form
\begin{eqnarray}
	\label{eq:eta_n}
	\left[\begin{array}{c} \eta_1 \\ \vdots \\ \eta_{n-1} \\ \eta_n \end{array}\right]=\left[\begin{array}{c} \sum_{j=1}^{n-1}\sum_{k=j-1}^{N}c_{jk}{k\choose j-1}\beta^k\lambda^{k+1-j}+\sum_{k=n-1}^{N}[\sum_{j=1}^{n-1}a_jc_{jk}]{k\choose n-1}\beta^k\lambda^{k+1-n} \\
	\vdots \\
	\sum_{k=0}^{N}c_{(n-1)k}{k\choose 0}\beta^k\lambda^{k}+\sum_{k=1}^{N}[\sum_{j=1}^{n-1}a_jc_{jk}]{k\choose 1}\beta^k\lambda^{k-1} \\
\sum_{k=0}^{N}[\sum_{j=1}^{n-1}a_jc_{jk}]\beta^k\lambda^{k} \end{array}\right].
\end{eqnarray}
Observe that $\sum_{k=0}^{N}c_{(n-1)k}\beta^k\lambda^k=\eta_{n-1}-\frac{d}{d\lambda}\eta_n$ and furthermore $\sum_{k=0}^{N}c_{(n-2)k}\b^k\lambda^k=\eta_{n-2}-\frac{d}{d\lambda}(\eta_{n-1}-\frac{d}{d\lambda}\eta_n)-\frac{1}{2!}\frac{d^2}{d\lambda^2}\eta_n$. Thus, for each $j=2,\ldots,n-1$, we can express $\sum_{k=0}^{N}c_{jk}\beta^k\lambda^k$ as a linear combination of $\eta_{j},\eta_{j+1},\ldots,\eta_n$ 
and their derivatives. In addition, from the last equation in \eqref{eq:eta_n}, we can express $\sum_{k=0}^{N}c_{1k}\beta^k\lambda^k=\frac{1}{a_1}(\eta_n-\sum_{j=2}^{n-1}a_j\sum_{k=0}^{N}c_{jk}\beta^k\lambda^k)$. As a result, $\sum_{j=1}^{n-1}\sum_{k=j-1}^{N}c_{jk}{k\choose j-1}\beta^k\lambda^{k+1-j}$ can be represented in terms of $\eta_2,\ldots,\eta_n$ and their derivatives. Together with the expression $\sum_{k=n-1}^{N}[\sum_{j=1}^{n-1}a_jc_{jk}]{k\choose n-1}\b^k\lambda^{k+1-n}=\frac{1}{(n-1)!}\frac{d^{n-1}}{d\lambda^{n-1}}\eta_n$, we conclude that $\eta_1$ is correlated to $\eta_2,\ldots,\eta_n$ and their higher order derivatives. Therefore, the system \eqref{eq:jordanprop} is not ensemble controllable, which implies that $\text{rank}(B)=n$ is necessary.
\hfill$\Box$

\begin{theorem}
	\label{thm:Jordan}
	Given the time-invariant linear ensemble system
	\begin{eqnarray*}
		\Sigma_2:\left\{
		\begin{array}{ll}
			\frac{d}{dt}X(t,\beta)=\beta AX(t,\beta)+BU(t), \\
			\beta\in K=[\b_1,\b_2]\subset\mathbb{R}^+ \ (\text{or}\ \mathbb{R}^-),
		\end{array}
		\right.
	\end{eqnarray*}
	where $X\in M\subset\mathbb{R}^n$, 
	the control $U:[0,T]\rightarrow\mathbb{R}^m$ is piecewise continuous, and $A\in\mathbb{R}^{n\times n}$ and $B\in\mathbb{R}^{n\times m}$ are constant matrices with $\text{rank}(A)=n$ and $\text{rank}(B)=m<n$. Suppose that $A$ is not diagonalizable and is similar to a Jordan form $J$ such that $A=PJP^{-1}$, where $J=\text{diag}(J_{1},...,J_{k})$, in which $J_i=\text{diag}(J_{i1},\ldots,J_{ir(i)})$, $i=1,\ldots,k$, and $J_{ij}\in\mathbb{R}^{d_{ij}\times d_{ij}}$, $j=1,\ldots,r(i)$, denotes a Jordan block with the eigenvalue $\lambda_i$. Note that $J_{ij}$ can be a scalar. Consider the transformed system $Y(t,\b)=P^{-1}X(t,\b)$, satisfying
\begin{eqnarray}
	\label{eq:sysjordan}
	\frac{d}{dt}Y(t,\b)=\b JY(t,\b)+\widetilde{B}U(t),
\end{eqnarray}
where $\widetilde{B}=P^{-1}B$. Writing $\widetilde{B}$ as the concatenation of matrices $\widetilde{B}_i$, $i=1,\ldots,k$, such that
$$\widetilde{B}=\left[\begin{array}{c}\widetilde{B}_1 \\ \hline \vdots \\ \hline\widetilde{B}_k\end{array}\right]\quad\text{and}\quad\widetilde{B}_i=\left[\begin{array}{c}\widetilde{B}_{i1} \\ \hline \vdots \\ \hline\widetilde{B}_{ir(i)}\end{array}\right]$$
where $\widetilde{B}_{ij}\in\mathbb{R}^{d_{ij}\times m}$ for $j=1,\ldots,r(i)$, the system $X$ is uniformly ensemble controllable if and only if
\begin{enumerate}
 		\item[(i)] $\text{rank}(\widetilde{B}_{ij})=d_{ij}$ for $i=1,\ldots,k$ and $j=1,\ldots,r(i)$;
		\item[(ii)] the corresponding diagonal system $\frac{d}{dt}Y(t,\b)=\b\Lambda Y(t,\b)+\widetilde{B}U(t)$ is uniformly ensemble controllable, where $\Lambda=\text{diag}(\Lambda_1,\ldots,\Lambda_k)$ in which $\Lambda_i=\text{diag}(\Lambda_{i1},\ldots,\Lambda_{ir(i)})$, $i=1,\ldots,k$, with $\Lambda_{ij}=\text{diag}(\lambda_i,\ldots,\lambda_i)\in\mathbb{R}^{d_{ij}\times d_{ij}}$ for $j=1,\ldots,r(i)$.
\end{enumerate}
\end{theorem}
{\it Proof:} \noindent (Sufficiency) Suppose that the conditions (i) and (ii) hold. Then, the rank condition in (i) implies that $\text{span}\{\b^s J_{ij}^s\, \tilde{b}_{ij}^{(l)}\}=\text{span}\{\b^s \Lambda_{ij}^s\, \tilde{b}_{ij}^{(l)}\}$, for all $i=1,\ldots,k$ and $j=1,\ldots,r(i)$, where $s=1,2,\ldots$ and $\tilde{b}_{ij}^{(l)}$ is the $l^{th}$ column of $\widetilde{B}_{ij}$ with $l=1,\ldots,m$. 
Therefore, $\overline{\{\b JY,\tilde{b}_j\}}_{LA}=\overline{\{\b \Lambda Y,\tilde{b}_j\}}_{LA}$, in which $\tilde{b}_j$ is the $j^{th}$ column of $\widetilde{B}$. This together with the condition (ii) guarantees uniform controllability of the system \eqref{eq:sysjordan}. 

(Necessary) Suppose that the system \eqref{eq:sysjordan} is ensemble controllable. This implies that each subsystem $\frac{d}{dt}{Y_{ij}}=\b J_{ij}Y_{ij}+\widetilde{B}_{ij}U$, $i=1,\ldots,k$ and $j=1,\ldots,r(i)$, is ensemble controllable, where $Y_{ij}\in\mathbb{R}^{d_{ij}}$, and hence $\text{rank}(\widetilde{B}_{ij})=d_{ij}$ according to Proposition \ref{prop:jordan}. Consequently,  
$\overline{\{\b JY,\tilde{b}_j\}}_{LA}=\overline{\{\b \Lambda Y,\tilde{b}_j\}}_{LA}$, $j=1,\ldots,m$, which results in (ii).\hfill$\Box$

We create the following example to illustrate the results of Theorem \ref{thm:Jordan}.

\begin{example} 	
	\label{ex:thm3}
	Consider the ensemble system $\frac{d}{dt}Y(t,\beta)=\beta JY(t,\beta)+\widetilde{B}U(t)$ with
	$$J=\left[\begin{array}{cccc} 1 & 0 & 0 & 0 \\ 0 & \alpha & 0 & 0 \\ 0 & 0 & 2 & 1 \\ 0 & 0 & 0 & 2 \end{array}\right],\quad \widetilde{B}=\left[\begin{array}{cc} 1 & 0 \\ 0 & 1 \\ 1 & 2 \\ 1 & 0 \end{array}\right],$$
and $\b\in K=[1,1.5]$. We consider two values of $\alpha$ to show different scenarios of ensemble controllability.\\
(Case I: $\alpha=2$) In this case, $J$ has two distinct eigenvalues, and then we let
$$\widetilde{B}_1=\left[\begin{array}{cc} 1 & 0\end{array}\right], \quad \widetilde{B}_{21}=\left[\begin{array}{cc} 0 & 1\end{array}\right], \quad \widetilde{B}_{22}=\left[\begin{array}{cc} 1 & 2 \\ 1 & 0 \end{array}\right].$$
It is clear that the condition (i) in Theorem \ref{thm:Jordan} is satisfied. However, the system $\frac{d}{dt}Y(t,\b)=\b\Lambda Y(t,\b)+\widetilde{B}U(t)$, where $\Lambda=\text{diag}(1,2,2,2)$, is not ensemble controllable according to Theorem \ref{thm:diagonal}. The violation of the condition (ii) in Theorem \ref{thm:Jordan} leads to the failure of ensemble controllability. \vspace{4pt}

\noindent (Case II: $\alpha=4$) In this case, it is easy to verify that the condition (i) holds. The corresponding diagonal system $\frac{d}{dt}Y(t,\b)=\b\Lambda Y(t,\b)+\widetilde{B}U(t)$ with $\Lambda=\text{diag}(1,4,2,2)$ is ensemble controllable according to Theorem \ref{thm:diagonal}. Thus, the system is ensemble controllable.
\end{example}

\subsubsection{$A$ Has Complex Eigenvalues}
Recall the time-invariant linear ensemble system $\frac{d}{dt}X(t,\b)=\b AY(t,\b)+BU(t)$ with $\b\in K\subset\mathbb{R}^+$ (or $\mathbb{R}^-$) as described in \eqref{eq:LTI_1}. Here, we consider the case when $A$ is full rank and, for simplicity, assume that $A$ is similar to a diagonal matrix $J$, i.e., $A=PJP^{-1}$, where $J=\text{diag}(\lambda_1,\ldots,\lambda_n)$, $\lambda_i\in\mathbb{C}$ for some $i=1,\ldots,n$, and $\lambda_i\neq 0$ for all $i=1,\ldots,n$. Without loss of generality, we may assume the first $2k$ eigenvalues are complex, namely, $\lambda_{2j-1}=\lambda_{2j}^{\dagger}$ for $j=1,\ldots,k$ with the corresponding eigenvectors $v_{2j-1}=v_{2j}^{\dagger}$, where $\dagger$ denotes the complex conjugate. As discussed in earlier sections, ensemble controllability is determined by whether any target ensemble state can be approximately reached. To fix the idea, we begin with analyzing a single-input two-dimensional system. 

Suppose that $A\in\mathbb{R}^{2\times 2}$ is similar to $J$, i.e., $A=PJP^{-1}$ where $J=\text{diag}(\lambda,\lambda^{\dagger})$ and $\lambda\in\mathbb{C}$. Then, the transformed system $Y(t,\b)=P^{-1}X(t,\b)$ follows $\frac{d}{dt}Y(t,\b)=\b JY(t,\b)+\tilde{b}u$, where $\tilde{b}=P^{-1}B=(\zeta,\zeta^{\dagger})'$, $\zeta\in\mathbb{C}$, and $\zeta\neq 0$. Then, the reachable set is characterized by
$\mathcal{L}_0=\text{span}\{c_k(\b\lambda)^k \zeta,\ c_k\in\mathbb{R},\ k=0,1,\ldots\}$. The system is ensemble controllable if, for any $\e>0$, there exist coefficients $c_k$ such that
\begin{eqnarray}
	\label{eq:approx1}
	\|\xi(\b)-\sum_{k=0}^{\infty} c_k(\b\lambda)^k\zeta\|_{\infty}<\e,
\end{eqnarray}
where $(\xi(\b),\xi^{\dagger}(\b))'=e^{-\b JT}Y_F(\b)-Y_0(\b)$, and $Y_0=P^{-1}X_0$ and $Y_F=P^{-1}X_F$ are the initial and the target state, respectively. Representing $\lambda=re^{i\theta}$, where $r>0$, the approximation for $\xi(\b)\in\mathbb{C}$ in \eqref{eq:approx1} is reduced to 
$$\Big\|\left[\begin{array}{c} f_1(\b)  \\ f_2(\b) \end{array}\right]-\sum_{k=0}^{\infty} c_k\left[\begin{array}{c} (\b r)^k\cos(k\theta)  \\ (\b r)^k\sin(k\theta) \end{array}\right]\Big\|_{\infty}<\e,$$
by defining $\frac{\xi(\b)}{\zeta}\doteq f_1(\b)+if_2(\b)$. The system is uniformly ensemble controllable if such polynomial approximation is possible for any given continuous functions $f_1,f_2:K\rightarrow\mathbb{R}$. For some special cases, e.g., $\theta=\frac{\pi}{2}$, or more generally $\theta=q\pi$ where $q$ is a rational number, it is evident that such a simultaneous polynomial approximation task using the common $c_k$ is possible. Because, for example, when $\theta=\frac{\pi}{2}$, i.e., $\lambda=ir$, the approximation to $f_1$ and to $f_2$ becomes independent, i.e., $f_1(\b)\approx c_0+c_2 r^2 \b^2+c_4 r^4 \b^4+\ldots$ and $f_2(\b)\approx c_1 r \b-c_3 r^3 \b^3+c_5 r^5 \b^5-\ldots$, for which there exist coefficients $c_k$ such that both approximations can be achieved.

The same logic can be applied to investigate higher dimensional cases, and controllability is then determined based on the same procedure as to evaluate the feasibility of the reduced polynomial approximation problems. Note that the number of controls corresponds to that of the sets of coefficients. For example, in the above two-dimensional case, there is only one set of coefficients $c_k$ that can be chosen for approximation since the system is of single input. 
This indicates a general situation where an ensemble system (with complex eigenvalues) is more likely controllable if the number of the controls, $m$, is more comparable to the dimension of the system, $n$, when $m<n$. Finally, we note that similar analysis can be carried out when the system matrix $A$ is similar to a Jordan form with complex eigenvalues.

\section{Ensemble Control Synthesis and Simulations}
\label{sec:example}
In this section, we present several practical ensemble systems and construct optimal ensemble controls using an SVD-based computational method developed in our previous work \cite{Li_ACC12SVD}.

\begin{figure}[t]
        \centering
        \subfigure[]{\includegraphics[scale=0.5]{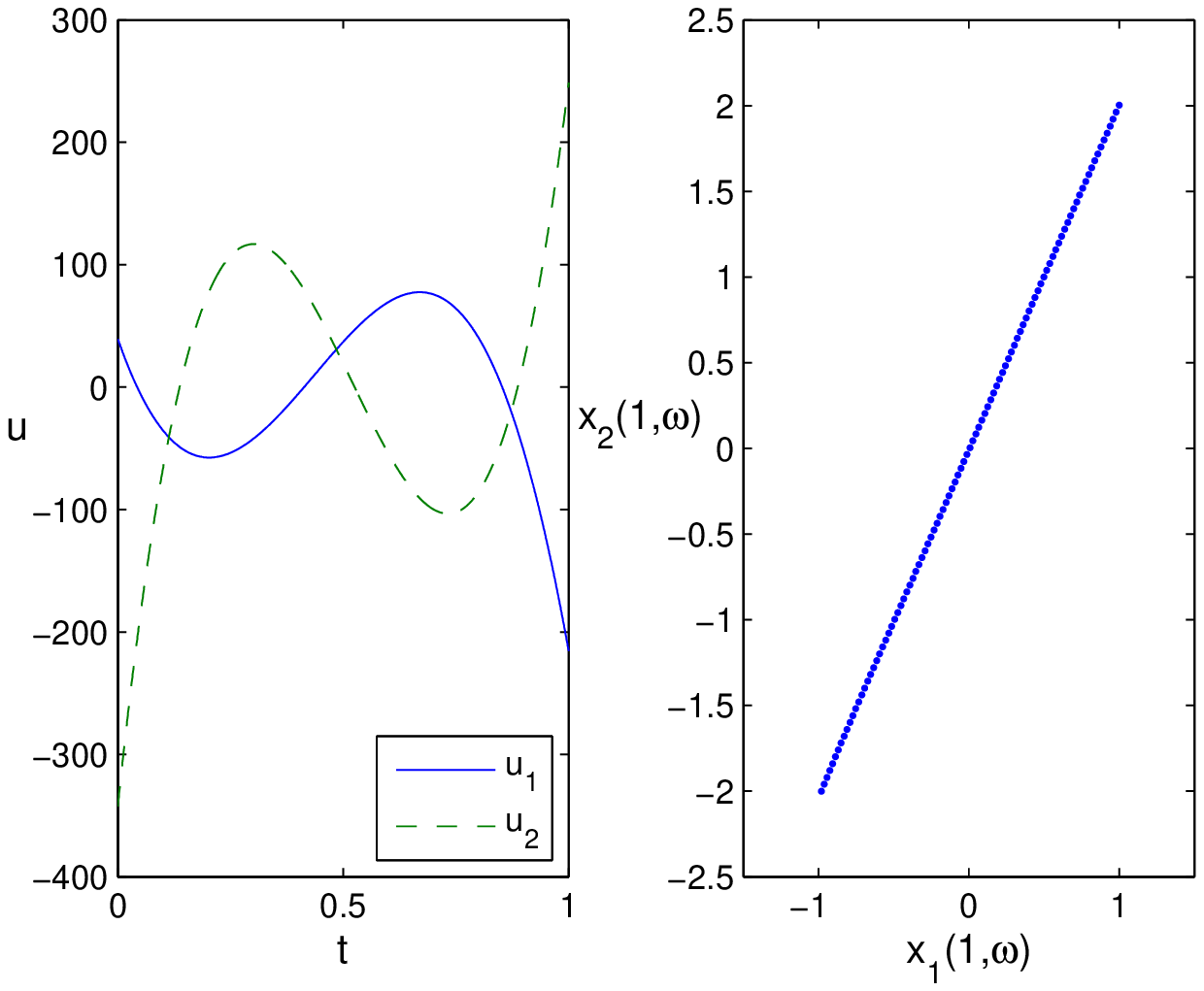}\label{fig:opt_control}}
        \subfigure[]{\includegraphics[scale=0.53]{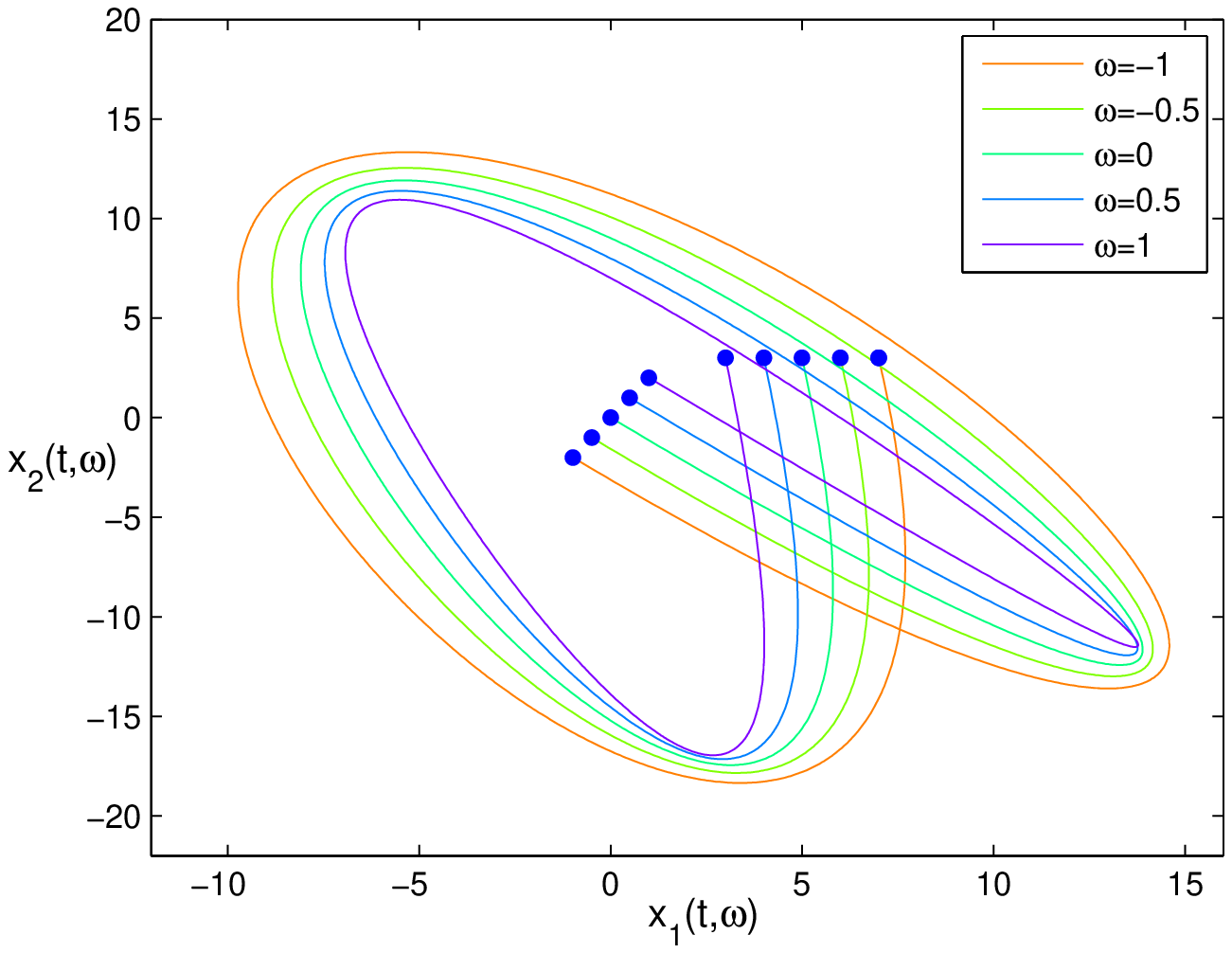}\label{fig:ho_traj}}
        \caption{\subref{fig:opt_control} The minimum-energy ensemble control law that steers the ensemble from $(5-2\w,3)'$ to $(\w,2\w)'$ at $T=1$ and the resulting final states $X(1,\w)$ for $\omega\in[-1,1]$. \subref{fig:ho_traj} Sample trajectories for $\omega=-1$, $-0.5$, $0$, $0.5$, and $1$.}
        \label{fig_har_osci}
\end{figure}

\begin{example}[Harmonic Oscillators]
	Consider steering an ensemble of harmonic oscillators modeled in \eqref{eq:harmonic} with their frequencies $\w\in K=[-1,1]$ from the initial state $X_0(\omega)=X(0,\w)=(5-2\w,3)'$ to the target state $X_F(\omega)=(\omega, 2\omega)'$ at time $T=1$. This system is uniformly ensemble controllable as shown in Example \ref{ex:harmonic}. The minimum-energy ensemble control law for steering this system and the resulting final states and sample trajectories are illustrated in Figure \ref{fig_har_osci}.
\end{example}

\begin{figure}[t]
        \centering
        \subfigure[]{\includegraphics[scale=0.55]{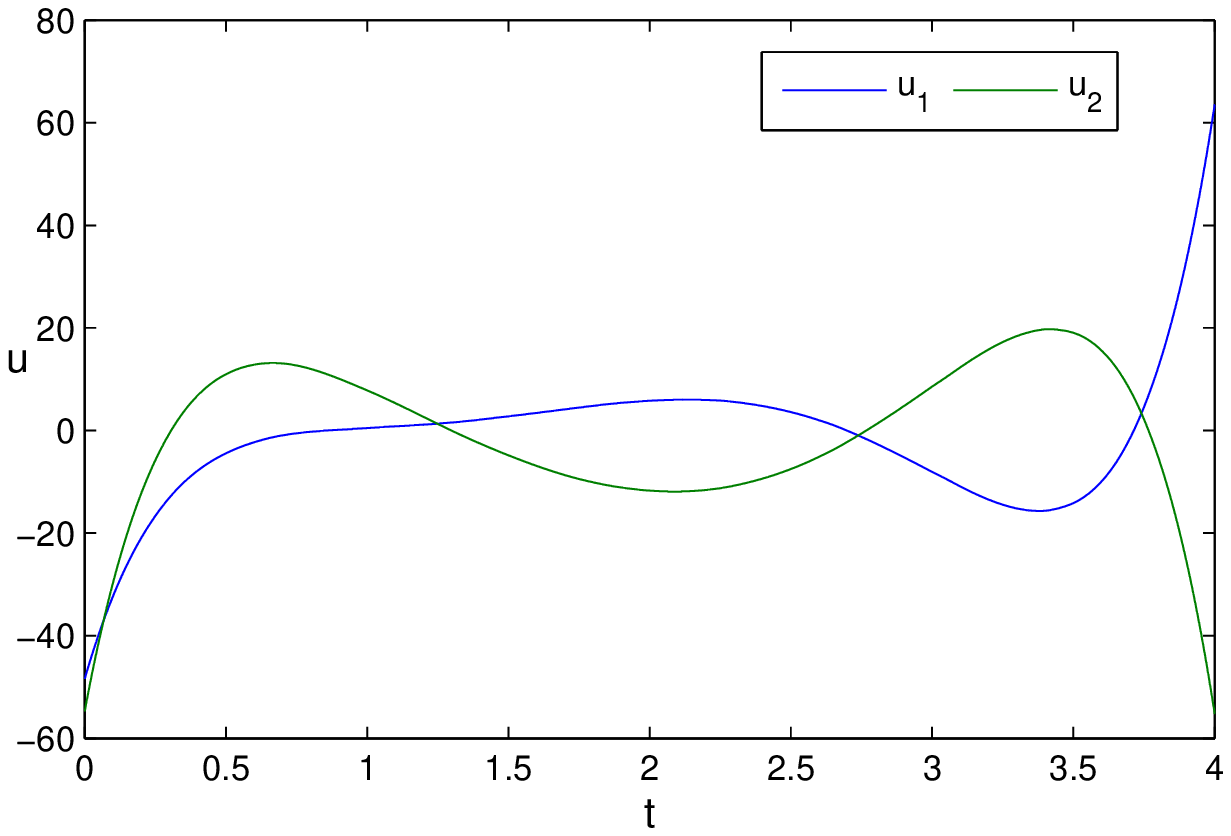}\label{fig:air_control2}}
        \subfigure[]{\includegraphics[scale=0.55]{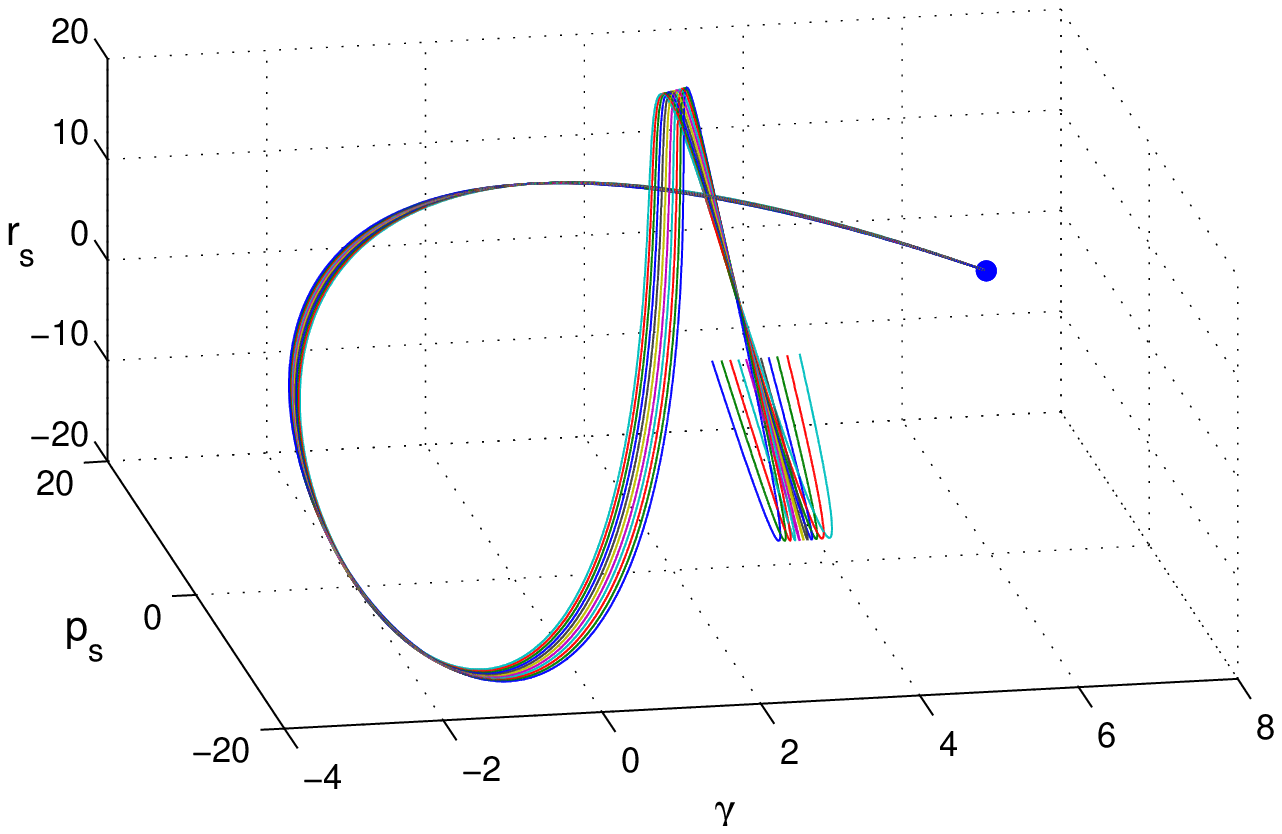}\label{fig:air_traj2}}
\caption{\subref{fig:air_control2} The minimum-energy ensemble control law that drives the aircraft system in the presence of uncertainty from $(2\pi,6,4)'$ to $(\pi\epsilon,\epsilon,0)'$ at $T=4$. \subref{fig:air_traj2} Sample trajectories for $\epsilon\in[0.8,1.2]$ following the optimal control shown in (a).}
        \label{fig:aircraft_II}
\end{figure}

\begin{example}[Aircraft System]
	\rm Consider the aircraft system described in Example \ref{ex:aircraft} with
$$A(\epsilon)=\epsilon A=\epsilon\left[\begin{array}{ccc} 0 & 0.5 & 1 \\ 2 & 0.5 & 0 \\ 0.5 & 0 & 0.5 \end{array}\right], \quad B=\left[\begin{array}{cc} 1 &  0 \\ 0 & 1 \\ 0 & 1\end{array}\right],$$
where $\epsilon\in K=[0.8,1.2]$. It is straightforward to verify that this system is uniformly ensemble controllable by using Theorem \ref{thm:diagonal}. The minimum-energy control that steers this ensemble system from $X_0=(2\pi,6,4)'$ to $X_F(\epsilon)=(\pi\epsilon,\epsilon,0)'$ and the resulting sample trajectories are shown in Figure \ref{fig:aircraft_II}.
\end{example}

\begin{figure}[t]
        \centering
        \subfigure[]{\includegraphics[scale=0.52]{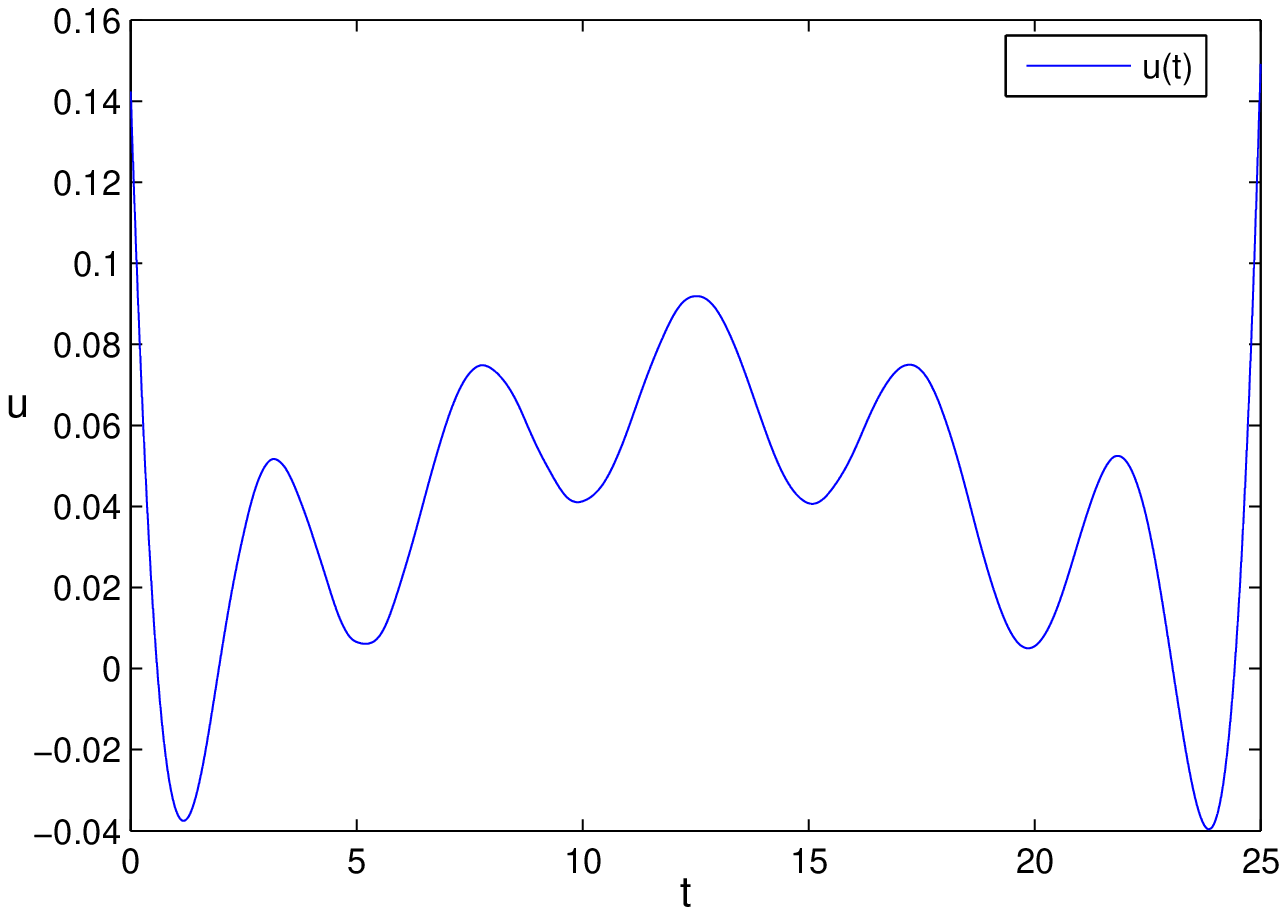}\label{fig:qtu1}}
        \subfigure[]{\includegraphics[scale=0.55]{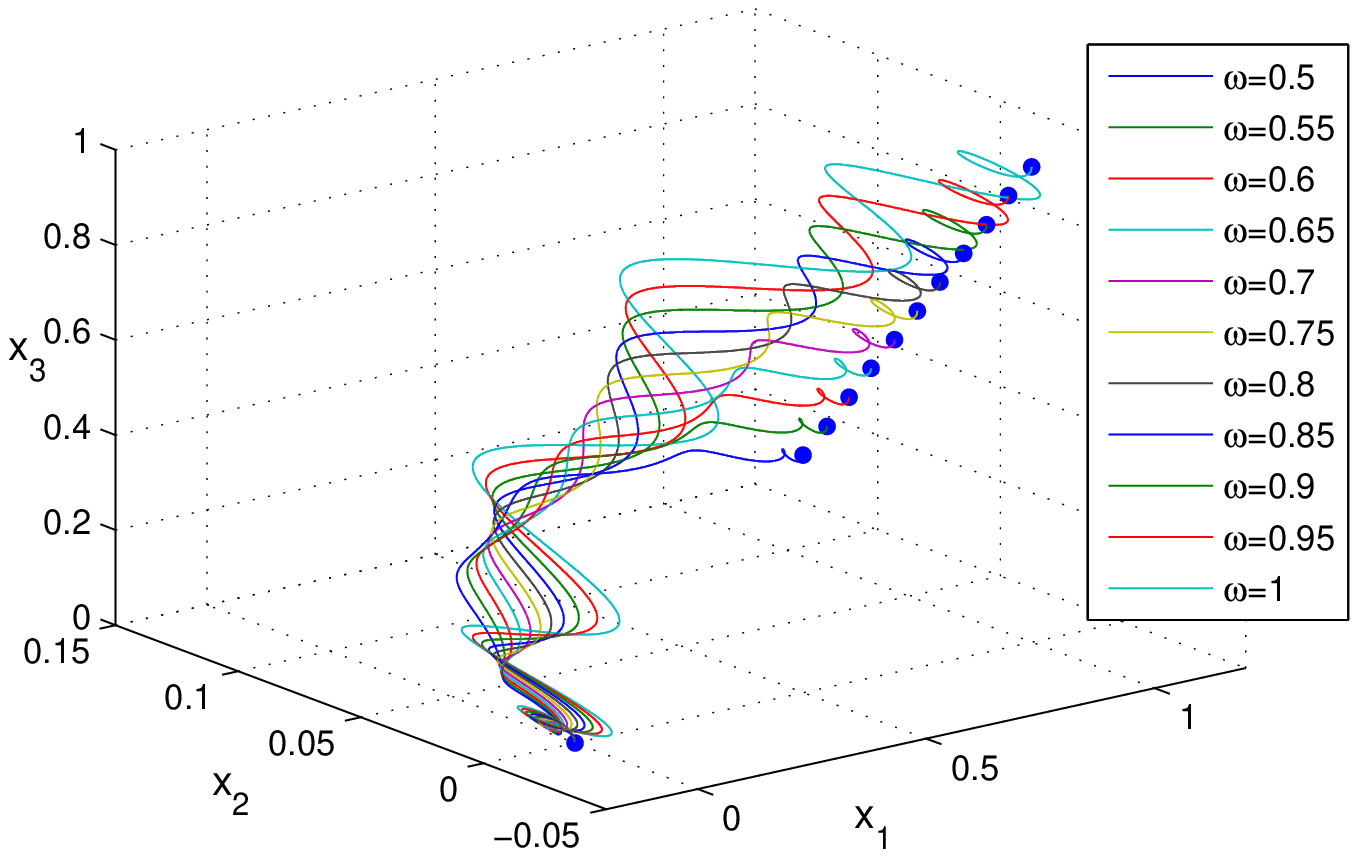}\label{fig:qt3}}
        \caption{\subref{fig:qtu1} The minimum-energy ensemble control that transports the atoms from $(0,0,0)'$ to $(\w,0,\w)'$ at $T=25$. 
\subref{fig:qt3} Sample trajectories for $\omega\in[0.5,1]$ following the optimal ensemble control law.}
        \label{fig:quantum}
\end{figure}

\begin{example}[Quantum Transport] 
	\label{ex:quantrans}
	Efficient transport of ultracold atoms trapped in a harmonic potential is an important goal in atomic physics, which leads to broad applications to basic sciences, metrology, and quantum information processing \cite{Murphy09, Hansel01}. The quantum mechanical description of frictionless transport of atoms can be reduced to the steering of a three-dimensional time-invariant linear system \cite{Li_TAC14_Transport}, given by $\frac{d}{dt}X(t,\w)=A(\w)X(t,\w)+B(\w)u(t)$ with the respective initial and terminal state, $X(0,\w)=(0,0,0)'$ and $X_F(\w)=(\w,0,\w)'$, where	
	$$A(\omega)=\left[\begin{array}{ccc} 0 & \omega & 0 \\ -\omega & 0 & \omega  \\ 0 & 0 & 0 \end{array}\right], \quad B(\w)=\left[\begin{array}{c} 0 \\ 0 \\ \omega \end{array}\right],$$
	and $\omega\in[0.5,1]$ represents the uncertainty in the angular frequency of the harmonic trapping potential. 
	Because a linear variation in $\w$ appears in the control term of the third state $x_3$, i.e., $\dot{x}_3=\w u$, this system is not ensemble controllable. 
	However, the reachable subspace can be characterized by the generators in the Lie algebra,
$$\mathcal{L}_0=\text{span}\{\left[\begin{array}{c} 0 \\ 0 \\ \omega \end{array}\right],\left[\begin{array}{c} 0 \\ \omega^{2k} \\ 0 \end{array}\right],\left[\begin{array}{c} \omega^{2k+1} \\ 0 \\ 0   \end{array}\right], k=1,2,3,\ldots\}.$$
We observe that $\xi(\w)=e^{-A(\w)T}X_F(\w)-X(0,\w)\in\overline{\mathcal{L}_0}$, and hence the state $(\omega,0,\omega)'$ is ensemble accessible from $(0,0,0)'$. The minimum-energy control that transports 
the atoms at time $T=25$ and the resulting sample trajectories are displayed in Figure \ref{fig:quantum}.
\end{example}

\section{CONCLUSIONS}
In this paper, we investigated a class of time-invariant linear ensemble systems whose natural dynamics are linear in the system parameter. We derived explicit controllability rank conditions that are easy to be checked. Specifically, we studied the cases when the parameter set is across the origin and when the parameter values are strictly positive or negative. Examples and numerical simulations were provided to demonstrate these theoretical results, and optimal ensemble control laws for steering the ensemble systems in these examples were constructed using an SVD-based computational method. We plan to generalize the developed methods, which were based on the notion of polynomial approximation, to establish explicit controllability conditions for general time-invariant linear ensemble system as presented in \eqref{eq:LTI}.


\section{APPENDIX}

\subsection{Proof of Proposition \ref{prop:range}}
\label{appd:proof}
\subsubsection{($\overline{R(\tilde{L})}\subseteq\overline{\mathcal{L}_0}$)}
Suppose that $\xi\in R(\tilde{L})$. Then, there exists some $U\in PC^m[0,T]$ such that $\xi=\tilde{L} U$. Then, we have from \eqref{eq:L_1}
\begin{align}
	\xi(\b) &= \int_0^T\tilde{\Phi}(0,\s,\beta)B(\b)U(\s)d\s \nonumber\\
	&= \int_0^T\sum_{n=0}^{\infty}\frac{(-\s)^n}{n!}A^n(\beta)B(\b)U(\s)d\s \nonumber\\
	\label{eq:sum_integral_2}
	&= \sum_{n=0}^{\infty}A^n(\beta)B(\b)\left[\int_0^T\frac{(-\s)^n}{n!}U(\s)d\s\right], 
\end{align}
This implies that $\xi(\b)\in\text{span}\{A^n(\b)b_j(\b)\}$, $n=0,1,\ldots$ and $j=1,\ldots,m$, where $b_j(\b)$ is the $j^{th}$ column of $B(\b)$, and hence $\xi\in\mathcal{L}_0$. It follows that $\overline{R(\tilde{L})}\subseteq\overline{\mathcal{L}_0}$. Note that the exchange of the infinite sum and the integration in \eqref{eq:sum_integral_2} is valid because of the uniform convergence of the series expansion of $\tilde{\Phi}$ and the boundedness of $B(\b)U(t)$ for $t\in [0,T]$ and for each $\b\in K$.

\subsubsection{($\overline{\mathcal{L}_0}\subseteq\overline{R(\widetilde{L})}$)}
\label{sec:proof1}
Consider an element $\eta_0\in\mathcal{L}_0$ of the form $\eta_0(\b)=\a_{ik}A^k(\beta)b_i(\b)$ for some $i\in\{1,\ldots,m\}$ and $k\in\{0,1,2,\ldots\}$, where $\a_{ik}\in\mathbb{R}$. Because $A\in C^{n\times n}(K)$, the transition matrix of the system \eqref{eq:LTI}, $\tilde{\Phi}(t,0,\beta)=e^{A(\b)t}$, is a convergent function on $D=[0,T]\times K$. Let $S_N(0,t,\beta)$ be a partial sum of the expansion of the inverse of the transition matrix, $\tilde{\Phi}(0,t,\b)$, defined by
$$S_N(0,t,\beta)=\sum_{j=0}^N \frac{A^j(\b) (-t)^j}{j\,!}$$
with $N>k$ such that $\|\widetilde{\Phi}(0,t,\beta)-S_N(0,t,\beta)\|_{\infty}<\e$. We now construct polynomial control functions $u_i=\sum_{j=0}^N a_{ij}t^j$, $a_{ij}\in\mathbb{R}$, such that
\begin{align}
	\label{eq:u1}
	\int_0^T \frac{(-t)^j}{j!}u_i(t)dt &=\left\{\begin{array}{ll}\a_{ik}, \quad \text{if}\quad j=k, \\0, \quad\ \ \, \text{if}\quad j\neq k, \end{array}\right.
\end{align}
where $1\leq i\leq m$ and $0\leq j\leq N$. This yields
\begin{eqnarray}
	\label{eq:error3}
	\int_0^T S_N(0,t,\beta)b_iu_i(t)dt=\a_{ik}A^k(\beta)b_i, \quad 1\leq i\leq m.
\end{eqnarray}
Thus, for $U=(0,\ldots,u_i,\ldots,0)'\in C^m([0,T])$, we have
\begin{align}
	&\|\tilde{L}U-\eta_0\|_{\infty} = \Big\|\int_0^T\widetilde{\Phi}(0,t,\beta)B(\b)U(t)dt-\a_{ik}A^k(\beta)b_i\Big\|_{\infty} \nonumber\\
	& \leq\Big\|\int_0^T\Big(\widetilde{\Phi}(0,t,\beta)-S_N(0,t,\beta)\Big)B(\b)U(t)dt\Big\|_{\infty} \nonumber\\
	& \ +\Big\|\int_0^T S_N(0,t,\beta)B(\b)U(t)dt-\a_{ik}A^k(\beta)b_i\Big\|_{\infty} \nonumber \\
	\label{eq:error2}
	& \leq\e T M_1 M_2
\end{align}
by the convergence of $\widetilde{\Phi}(t,0,\b)$ and by \eqref{eq:error3}, 
where $M_1=\|U\|_{\infty}$ and $M_2=\|B\|_{\infty}$, and the inequalities above are pointwise. This implies that $\eta_0\in \overline{R(\widetilde{L})}$. The same procedure as described above can be applied to show that any element $\eta=\sum_{k=0}^{\infty}\sum_{i=1}^m \a_{ik}\b^kA^kb_i\in\mathcal{L}_0$ is also an element of $\overline{R(\widetilde{L})}$. Consequently, we conclude that $\overline{\mathcal{L}_0}\subseteq\overline{R(\widetilde{L})}$. \hfill$\Box$

\subsection{M\"{u}ntz$-$Sz\'{a}sz Theorem \cite{Muntz1914}}
\label{appd:Muntz}
Let $\{\lambda_i\}_{i=1}^{\infty}$ be a sequence with $\inf_i\lambda_i>0$. Then
$$\text{span}\{1,x^{\lambda_1},x^{\lambda_2},\ldots\},$$
is dense in $C[0,1]$ if and only if
$$\sum_{i=1}^{\infty}\frac{1}{\lambda_i}=\infty.$$

%

%
%

\bibliographystyle{ieeetr}
\footnotesize
\bibliography{linear_controllability}

\begin{thebibliography}{10}

\bibitem{Li_PRA06}
J.-S. Li and N.~Khaneja, ``Control of inhomogeneous quantum ensembles,'' {\em
  Physical Review A}, vol.~73, p.~030302, 2006.

\bibitem{Li_PNAS11}
J.-S. Li, J.~Ruths, T.-Y. Yu, H.~Arthanari, and G.~Wagner, ``Optimal pulse
  design in quantum control: A unified computational method,'' {\em Proceedings
  of the National Academy of Sciences}, vol.~108, no.~5, pp.~1879--1884, 2011.

\bibitem{Glaser98}
S.~J. Glaser, T.~Schulte-{Herbr\"uggen}, M.~Sieveking, O.~Schedletzky, N.~C.
  Nielsen, O.~W. {S{\o}rensen}, and C.~Griesinger, ``Unitary control in quantum
  ensembles, maximizing signal intensity in coherent spectroscopy,'' {\em
  Science}, vol.~280, pp.~421--424, 1998.

\bibitem{Benabid10}
A.~L. Benabid, P.~Pollak, D.~Haffmann, C.~Gervason, M.~Hommek, J.~Perret,
  J.~de~Rougemont, and D.~Gao, ``Long-term suppression of tremor by chronic
  stimulation of the ventral intermediate thalamic nucleus,'' {\em The Lancet},
  vol.~337, pp.~403--406, 1991.

\bibitem{Foutz10}
T.~J. Foutz and C.~C. McIntyre, ``Evaluation of novel stimulus waveforms for
  deep brain stimulation,'' {\em Journal of Neural Engineering}, vol.~7, no.~6,
  p.~066008, 2010.

\bibitem{Li_JNE12}
A.~Zlotnik and J.-S. Li, ``Optimal entrainment of neural oscillator
  ensembles,'' {\em Journal of Neural Engineering}, vol.~9, no.~4, p.~046015,
  2012.

\bibitem{Becker12rob}
A.~Becker and T.~Bretl, ``Approximate steering of a unicycle under bounded
  model perturbation using ensemble control,'' {\em IEEE Transactions on
  Robotics}, vol.~28, no.~3, pp.~580--591, 2012.

\bibitem{Li_TAC09}
J.-S. Li and N.~Khaneja, ``Ensemble control of bloch equations,'' {\em IEEE
  Transactions on Automatic Control}, vol.~54, pp.~528--536, 2009.

\bibitem{Li_TAC11}
J.-S. Li, ``Ensemble control of finite-dimensional time-varying linear
  systems,'' {\em IEEE Transastions on Automatic Control}, vol.~56,
  pp.~345--357, 2011.

\bibitem{Li_TAC13}
J.-S. Li, I.~Dasanayake, and J.~Ruths, ``Control and synchronization of neuron
  ensembles,'' {\em IEEE Transactions on Automatic Control}, vol.~58, no.~8,
  pp.~1919--1930, 2013.

\bibitem{Li_NOLCOS10}
J.-S. Li, ``Control of a network of spiking neurons,'' in {\em 8th IFAC
  Symposium on Nonlinear Control Systems}, (Bologna, Italy), Sep. 2010.

\bibitem{Li_PRL13}
A.~Zlotnik, Y.~Chen, I.~Z. Kiss, H.-A. Tanaka, and J.-S. Li, ``Optimal waveform
  for fast entrainment of weakly forced nonlinear oscillators,'' {\em Physical
  Review Letters}, vol.~111, no.~2, p.~024102, 2013.

\bibitem{Wilson14}
D.~Wilson and J.~Moehlis, ``Optimal chaotic desynchronization for neural
  populations,'' {\em {SIAM Journal on Applied Dynamical Systems}}, vol.~13,
  no.~1, pp.~276--305, 2014.

\bibitem{Li_TAC12_QCP}
J.~Ruths and J.-S. Li, ``Optimal control of inhomogeneous ensembles,'' {\em
  IEEE Transactions on Automatic Control: Special Issue on Control of Quantum
  Mechanical Systems}, vol.~57, no.~8, pp.~2021--2032, 2012.

\bibitem{Li_ACC12SVD}
A.~Zlotnik and J.-S. Li, ``Synthesis of optimal ensemble controls for linear
  systems using the singular value decomposition,'' in {\em American Control
  Conference}, (Montreal, Canada), pp.~5849--5854, Jun. 2012.

\bibitem{Li_TAC14_Transport}
D.~Stefanatos and J.-S. Li, ``Minimum-time quantum transport with bounded trap
  velocity,'' {\em IEEE Transactions on Automatic Control}, vol.~59, no.~3,
  pp.~733--738, 2014.

\bibitem{Brockett10flock}
R.~Brockett, ``On the control of a flock by a leader,'' {\em {Proceedings of
  the Steklov Institute of Mathematics}}, vol.~268, no.~1, pp.~49--57, 2010.

\bibitem{Liu11}
Y.-Y. Liu, J.-J. Slotine, and A.-L. Barabasi, ``Controllability of complex
  networks,'' {\em {Nature}}, vol.~473, no.~1, pp.~167--173, 2011.

\bibitem{Hochberg06}
L.~R.~H. et~al., ``Neuronal ensemble control of prosthetic devices by a human
  with tetraplegia,'' {\em Nature}, vol.~442, pp.~164--171, 2006.

\bibitem{Beauchard10}
K.~Beauchard, J.-M. Coron, and P.~Rouchon, ``Controllability issues for
  continuous-spectrum systems and ensemble controllability of bloch
  equations,'' {\em Communications in Mathematical Physics}, vol.~296, no.~2,
  pp.~525--557, 2010.

\bibitem{Gohberg03}
I.~Gohberg, S.~Goldberg, and M.~A. Kaashoek, {\em Basic Classes of Linear
  Operators}.
\newblock Boston, MA: Birkh\"auser Verlag, 2003.

\bibitem{Lavretsky12}
E.~Lavretsky and K.~Wise, {\em Robust and Adaptive Control: With Aerospace
  Applications (Advanced Textbooks in Control and Signal Processing)}.
\newblock Springer-Verlag, 2012.

\bibitem{Murphy09}
M.~Murphy, L.~Jiang, N.~Khaneja, and T.~Calarco, ``High-fidelity fast quantum
  transport with imperfect controls,'' {\em Physical Review A}, vol.~79,
  p.~020301, Feb. 2009.

\bibitem{Hansel01}
W.~H\"ansel, J.~Reichel, P.~Hommelhoff, and T.~W. H\"ansch, ``Magnetic conveyor
  belt for transporting and merging trapped atom clouds,'' {\em Physical Review
  Letters}, vol.~86, pp.~608--611, Jan. 2001.

\bibitem{Muntz1914}
C.~H. M{\"{u}}ntz, ``{\"{U}}ber den approximationssatz von weierstrass,'' {\em
  H. A. Schwarz's Festschrift}, pp.~303--312, 1914.

\end{thebibliography}

\end{document}